\PassOptionsToPackage{unicode}{hyperref}
\PassOptionsToPackage{naturalnames}{hyperref}

\documentclass[a4paper,11pt]{article}

\usepackage{amsfonts,amssymb,amsmath,amsthm}	
\usepackage[english]{babel}                 	
\usepackage[T1]{fontenc}                    	
\usepackage{setspace}      						
\usepackage{graphicx}      						
\usepackage{tikz}          						
\usepackage{enumitem}
\usepackage[permil]{overpic}       				
\usepackage{fancybox}							
\usepackage{ifthen}								
\usepackage{subcaption}
\usepackage[ plainpages = false, pdfpagelabels,
	     bookmarks,
	     bookmarksopen = true,
	     bookmarksnumbered = true,
	     breaklinks = true,
	     linktocpage,
	     pagebackref,         
	     colorlinks = true,  
         hypertexnames=false,   
	     linkcolor = green!50!black,
	     urlcolor  = blue!50!black,
	     citecolor = blue!50!black,
	     anchorcolor = blue,
	     hyperindex = true,
	     hyperfigures
	     ]{hyperref}
\usepackage{tikz}
\usepackage{tikz-3dplot}
\usepackage[numbers]{natbib}

\usepackage{bibentry} 
\usepackage{fancyhdr, graphicx}
\usepackage{pgfgantt}
\usepackage{bm}

\usepackage{bbm}
\usepackage[colorinlistoftodos,prependcaption]{todonotes}
\usepackage[capitalise]{cleveref}
\usepackage{autonum} 
\usepackage{amssymb}
\usepackage{titlesec}
\usepackage{booktabs}
\usepackage[parfill]{parskip} 

\usepackage{pgfplots}
\pgfplotsset{compat=newest}
\usetikzlibrary{plotmarks}
\usetikzlibrary{arrows.meta}
\usepgfplotslibrary{patchplots}
\usepackage{grffile}
\usepackage{amsmath}
\pgfplotsset{plot coordinates/math parser=false}
\newlength\figureheight
\newlength\figurewidth
\pgfplotsset{every axis/.append style={
                    label style={font=\small},
                    tick label style={font=\footnotesize},
                    legend style={font=\tiny},
                    title style={font=\small}
                    }}
\usetikzlibrary{plotmarks}
\pgfplotsset{minor grid style={dotted,gray}} 
\pgfplotsset{every axis/.append style={thick, tick style=semithick}}
\pgfplotsset{every mark/.append style={mark size=1pt}}
\usepgfplotslibrary{groupplots}
\pgfplotsset{xticklabel style={/pgf/number format/fixed,
        /pgf/number format/precision=5},yticklabel style={/pgf/number format/fixed,
        /pgf/number format/precision=5}}

\newcommand{\includetikz}[1]{%
        \includegraphics{#1} 
}

\renewcommand{\d}{\, \mathrm{d}}                                   
\newcommand{\vect}[1]{\textrm{\boldmath${#1}$}}                 
\newcommand{\abs}[1]{|{#1}|}                 
\newcommand{\norm}[1]{\Vert{#1}\Vert}

\newcommand{\energy}{\mathcal{E}}
\newcommand{\mass}{\mathcal{M}}
\newcommand{\density}{\rho}
\newcommand{\flux}{j}
\newcommand{\momentum}{\mathcal{P}}
\newcommand{\Etot}{\energy}                 
\newcommand{\Epot}{\energy_\mathrm{pot}}                 
\newcommand{\Ekin}{{\energy}_\mathrm{kin}}               
\renewcommand{\u}{\vect{u}}                              
\newcommand{\trec}{\tau_\mathrm{rec}}      

\newcommand{\map}{\vect{X}}
\newcommand{\Ncoarse}{N}
\newcommand{\Npsi}{N_{\psi}}
\newcommand{\Nfine}{N_\text{f}}
\newcommand{\vmap}{\map}
\newcommand{\xmap}{\vect{\Phi}}

\newcommand{\email}[1]{\href{mailto:#1}{\texttt{#1}}}           
\newcommand{\RomanNumeralCaps}[1]                               
    {\MakeUppercase{\romannumeral #1}}
\newcommand{\define}[1]{\textit{#1}}                            

{\theoremstyle{definition}\newtheorem{definition}{Definition}[section]}

{\theoremstyle{remark} }            

{\theoremstyle{lemma}
}

\newcommand{\rev}[1]{#1}

\graphicspath{{figures/two_stream/}{figures/landau_damping/}}

\binoppenalty=\maxdimen 
\relpenalty=\maxdimen   

\usepackage[margin=1in]{geometry}
\date{}

\title{\textbf{A Characteristic Mapping Method for Vlasov--Poisson with Extreme Resolution Properties}}
\author{
Philipp Krah\footnote{Institut de Mathématiques de Marseille, Aix-Marseille Université, \email{philipp.krah@univ-amu.fr}},
Xi-Yuan Yin\footnote{LMFA-CNRS, Ecole Centrale de Lyon, \email{xi-yuan.yin@ec-lyon.fr}},
Julius Bergmann\footnote{Institut de Mathématiques de Marseille, Aix-Marseille Université, \email{julius@bei-bergmann.de}},\\
Jean-Christophe Nave\footnote{McGill University, \email{jcnave@math.mcgill.ca}} \, and
Kai Schneider\footnote{Institut de Mathématiques de Marseille, Aix-Marseille Université, \email{kai.schneider@univ-amu.fr}}
}

\nobibliography{references} 
\begin{document}

\maketitle

\begin{abstract}
We propose an efficient semi-Lagrangian characteristic mapping method for solving the one+one-dimensional Vlasov--Poisson equations with high precision on a coarse grid. The flow map is evolved numerically and exponential resolution in linear time is obtained.
Global third-order convergence in space and time is shown and conservation properties are assessed. For benchmarking, we consider linear and nonlinear Landau damping and the two-stream instability. We compare the results with a Fourier pseudo-spectral method and results from the literature. The extreme fine-scale resolution features are illustrated showing the method's capabilities to efficiently treat filamentation in fusion plasma simulations.
\end{abstract}



\section{Introduction}

Flows in the general setting transport quantities, e.g. fluids, plasmas, particles or their probability distribution function (PDF), from one place to another. Typically they generate rich mathematical multi-scale structures even from simple analytical initial conditions and require well-adapted numerical methods to solve the underlying governing partial differential equations.
In the present work, we focus on the Vlasov--Poisson (VP) system, modeling particle evolution under their self-consistent electric field, and neglecting their collisions. We consider the equation in 1d+1d, i.e. in the space of positions and velocities. This transport equation in phase space has numerous applications in plasma and astrophysics \cite{Palmroth2018}. \rev{A} prominent example is Landau damping in statistical physics of ionized gases.

As a kinetic equation,
the VP equation is intrinsically high-dimensional, however, it is especially challenging due to the filamentation of the particle distribution function in phase space. The filamentation is a physical property of the system, that leads to a continuous generation of fine scales with increasing time.
Resolving these scales requires computer simulations to accumulate more and more computer storage and resources within time.  

To overcome the filamentation problem in numerical simulations two main approaches have been followed: the Lagrangian approach and the semi-Lagrangian approach. These methods shift the evolution of the PDF towards the evolution of its underlying transport structure.
In the purely Lagrangian approach, known as the particle-in-cell method (PIC) (see for review \cite{HockneyEastwood2021,BirdsallLangdon2018,Verboncoeur2005}), a set of $N_\text{p}$ particles is used to statistically approximate the distribution function in velocity space. The particles move according to their equations of motion in an electromagnetic field that is interpolated from a fixed grid to the particle's position. As the PDF is determined statistically the overall result is subject to statistical noise, with the variance decreasing only slowly with $1/\sqrt{N_\text{p}}$. 
Different approaches for denoising have been proposed, see e.g. \citep{picklo2023denoising}, among them are wavelet-based density estimation \citep{gassama2007wavelet,VanYenDelCastilloNegretteSchneiderFargeChen2010}.

In order to obtain even more accurate results, while simultaneously respecting the particle's equation of motion, semi-Lagrangian approaches have been developed. This special type of Eulerian method represents the PDF on a fixed grid whereas time is evolved with the help of particle trajectories, known as the characteristic curve. Along these curves, the PDF is conserved and thus constant in time. The PDF is computed by tracing back the trajectory and interpolating its origin back on an Eulerian mesh. The first approach in this direction was already presented in the seventies by Cheng and Knorr \cite{ChenKnorr1976} and has ever since got continuous attention in research \cite{SonnendrückerRocheBertrandGhizzo1999,MangeneyCalifanoCavazzoniTravnicek2002,BesseSonnendruecker2003,BesseMehrenberger2008,CrouseillesRespaudSonnendruecker2009,CrouseillesMehrenbergerSonnendruecer2010,RossmanithSeal2011,CrouseillesMehrenbergerVecil2011,EinkemmerOstermann2014}. 

In \cite{SonnendrückerRocheBertrandGhizzo1999} Sonnendrücker et~al. formalized the semi-Lagrangian approach for solving the VP equations by employing spline interpolations at the feet of the characteristics. While this method is renowned for its high precision, it is accompanied by the need to solve a global tri-diagonal system induced by the spline interpolation scheme. To overcome this computational challenge, \cite{BesseSonnendruecker2003} introduced local interpolation schemes based on Hermite and Lagrange polynomials, which were shown to exhibit high-order convergence results in space and second-order accuracy in time, as documented by \cite{BesseMehrenberger2008}.
The introduction of local interpolation schemes enhanced the parallelizability and adaptability of the method to unstructured meshes \cite{BesseSonnendruecker2003}. However, it introduced the requirement to transport gradients as additional fields. \rev{Subsequent developments in this direction focused on maintaining positivity and mass conservation, as presented for the backward \cite{CrouseillesMehrenbergerSonnendruecer2010}, and the forward semi-Lagrangian method \cite{CrouseillesRespaudSonnendruecker2009}.}

Another prominent methodology in the semi-Lagrangian setting is built on discontinuous Galerkin (DG) schemes, which have been explored by several research groups \cite{CrouseillesMehrenbergerVecil2011,EinkemmerOstermann2014,RossmanithSeal2011,QiuShu2011}. In DG-based semi-Lagrangian methods, the function is advected and subsequently projected onto the space of piecewise polynomial functions. DG methods offer the advantage of high-order local reconstruction and inherent mass conservation, as well as computational benefits when compared to cubic spline-based semi-Lagrangian methods, as reported by Einkemmer in \cite{Einkemmer2019}.

Although the scope of this paper lies on a semi-Lagrangian method we also briefly 
mention other methods like the purely Eulerian methods. They are based on a grid formulation of the distribution function and describe its time evolution on these fixed locations \cite{Klimas1987,KlimasFarrell1994,ZakiGardnerBoyd1988I,ZakiGardnerBoyd1988II,CaiLiWang2013,ParkerDellar2015,FilbetXiong2020}. In these works often the phase space is discretized with a tailored basis, or filtering is applied to handle filamentation. A comparison of semi-Lagrangian and Eulerian methods can be found in \cite{FilbetSonnendruecker2003}. A general overview of numerical techniques to solve VP systems is given in \cite{DimarcoPareschi2014}.

The characteristic mapping method (CMM) we are developing here for Vlasov--Poisson is a semi-Lagrangian method for solving transport-dominated problems. The origins of CMM were starting with semi-Lagrangian Gradient Augmented Level Set and Jet-Schemes \rev{\cite{NaveRosalesSeibold2010}}. Then an extension of the characteristic mapping method for linear advection equations was proposed in \cite{kohno2013new} and in \cite{mercier2020characteristic}.
One of the main features of CMM, which is a geometric approach and fully exploits the characteristic structure of the flow, is that
quantities advected by the flow can be constructed as the function composition of the initial condition with the characteristic map. This flow map is a diffeomorphism and for incompressible flows volume preserving. Moreover, it has a one-parameter semi-group structure (time) which allows decomposing the map into sub-maps. In contrast to other refinement methods we thus use compositional refinement. Therefore the PDF can be in principle evaluated at any given resolution at any given time, preserving the features of the initial condition like positivity. 
This however comes with the cost of storing more and more sub-maps, which may limit the total run-time.
More recently CMM has been extended for the
incompressible 2d Euler equations \cite{yin2021characteristic}, and the 3d Euler equations \cite{yin2023characteristic}. A generalization for solving flows on manifolds was likewise proposed for instance for computing tracer transport on the sphere \cite{taylor2023projection}
and for the incompressible Euler equations on the rotating sphere \cite{taylor2023characteristic}.

The current work aims to extend the CMM approach for the Vlasov--Poisson equation. 
We note that a related approach exists already in the literature \cite{KirchhartWilhelm2023}, which uses the compositional structure of a flow map similarly to solve Vlasov--Poisson. However, the details of the implementation differ as they are using the symplectic structure derived from a Hamilton formulation of the problem.

In this work, we perform some numerical analysis to assess the convergence order and show the third order in space and time. Code validation is provided to illustrate the theoretically shown convergence properties numerically. Benchmarking and comparison with a classical Fourier pseudo-spectral code yield cross-validation and show the efficiency of CMM with respect to its resolution capabilities. Finally, some high-resolution computations of the two-stream instability and Landau damping provide deep insight into the fine-scale structure of the distribution function.

The remainder of the paper is organized as follows. First, we recall the Vlasov--Poisson equation and summarize its key properties (\cref{sec:vp}). Then we describe the principles of the characteristic mapping method and its adaptation to Vlasov--Poisson, the numerical implementation and some numerical analysis in \cref{sec:CMM}.
\Cref{sec:pseudo-spectral} presents the pseudo-spectral discretization for Vlasov--Poisson used for comparison and benchmarking with CMM.
Different numerical test cases are presented in \cref{sec:numtest}.
Finally, some conclusion is given in \cref{sec:conclusion}.

\section{The Vlasov--Poisson Equation}
\label{sec:vp}
We introduce the Vlasov--Poisson equation 
\begin{align}
    \partial_t f + v \partial_x f + \partial_x \phi \partial_v f = 0\,,\label{eq:vlasov-advect}\\
    \partial_{xx} \phi +1 - l \int f \d{v} = 0\,, \label{eq:Vlasov--Poisson}
\end{align}
\rev{subject to an initial condition on $f$ and boundary conditions on $f$ and $\phi$}.\\
Here $f\colon\Omega_x\times \Omega_v\times[0,T]\to  \mathbb{R}$ denotes the \define{particle distribution function} (PDF), $\phi\colon\Omega_x\times[0,T]\to \mathbb{R}$ the \define{electric potential} with corresponding \define{electric field} $E(x,t) = -\partial_x \phi(x,t)$. 

The electron charge density of the plasma is $\density(x,t) = \int_{\Omega_v} f(x,v,t) \d v$ .
\begin{definition}[Mass and momentum]
For a solution $f(t,x,v)$ of the Vlasov--Poisson system we define the
\define{mass (0. moment)}
\begin{equation}
    \mass(t) =\int_{\Omega_x} \density(x,t) \d x = \iint_{\Omega_x\times \Omega_v}  f(x,v,t) \d x\, \d v\,,
\end{equation}
and its \define{moment (1. moment)}
\begin{equation}
    \momentum(t) =\int_{\Omega_x} \flux(x,t) \d x = \iint_{\Omega_x\times \Omega_v}  vf(x,v,t) \d x\, \d v\,.
\end{equation}
\end{definition}
For later reasons we write \cref{eq:vlasov-advect} in compact form using the Vlasov velocity field ${\u}=(v,\partial_x \phi)$ and the Vlasov nabla operator $\nabla^{\u}=(\partial_x,\partial_v)$. 
Thus we can rewrite \cref{eq:vlasov-advect} in terms of
\begin{equation}
\label{eq:vlasov-advectb}
    \partial_t f + {\u}\cdot \nabla^{\u} f = 0\,.
\end{equation}
    The velocity field ${\u}=(v,\partial_x \phi)$ 
    is divergence-free, which is straightforward to see:
    \begin{align}
    \nabla^{\u} \cdot {\u} = \partial_x v + \partial_v \partial_x \phi(x,t) = 0\,.
\end{align}

\begin{definition}[Energy]
For a solution $f(t,x,v)$ of the Vlasov--Poisson system we define its
\define{kinetic energy}
\begin{equation}
    \Ekin(t) =\frac{1}{2}\int_{\Omega_v} \int_{\Omega_x}  f(x,v,t) \abs{v}^2 \d x\, \d v\,,
\end{equation}
its \define{potential energy}
\begin{equation}
    \Epot(t) =\frac{1}{2}\int_{\Omega_x} \abs{E(x,t)}^2 \d x = \frac{1}{2}\int\int_{\Omega_x} \rho(t,x) \phi(t,x) \d{x}
\end{equation}
and \define{total energy}
\begin{equation}
    \energy(t) =\Epot(t)+ \Ekin(t)\,.
\end{equation}
\end{definition}

The total mass, momentum and total energy are conserved:
\begin{equation}
     \frac{\d \mass}{\d t}(t) = 0,\quad \frac{\d \momentum}{\d t}(t) = 0,\quad \frac{\d \energy}{\d t}(t) = 0\,.
\end{equation}

\section{Characteristic Mapping Method for Vlasov--Poisson}
\label{sec:CMM}
In the following, we describe how the characteristic mapping method can be used to solve the Vlasov--Poisson equation.
In \cref{subsec:CMM-basics} we review the fundamentals of the method for a scalar advection equation with a given velocity field \rev{\cite{kohno2013new,mercier2020characteristic}}. Thereafter \cref{subsec:advect-distr-func} introduces the velocity field and the boundary treatment that is specific to the Vlasov--Poisson equation.

\subsection{Characteristic Mapping Method}
\label{subsec:CMM-basics}
\newcommand{\scalar}{\theta}
\newcommand{\domain}{\Omega}

The CMM considers a non-linear advection equation
\begin{align}
\label{eq-def:CMM-advect}
    \begin{cases}
    \partial_t \scalar + {\u}\cdot \nabla\scalar =0 & \text{for } (\vect{x},t)\in \domain\times\mathbb{R}_+\\
    \scalar(\vect{x},0) = \scalar_0(\vect{x})           &  \vect{x}\in \domain\,,
    \end{cases}
\end{align}
with velocity field ${\u}\in \mathbb{R}^d$ that may depend on the advected state $\scalar\colon (\vect{x},t)\colon \domain\times\mathbb{R}_+\to \mathbb{R}$ itself.  
If we follow the points $(\vect{x},t)=(\vect{\gamma}(t),t)$ that satisfy the ODE:
\begin{align}
        \label{eq-def:characteristic_curve}
        \frac{\d \vect{\gamma}}{\d t}\, = {\u}   \quad \text{with} \quad
        \vect{\gamma}(0) = \vect{x}\in\domain
\end{align}
we have:
\begin{equation}
    \frac{\d {\scalar}}{\d t}(\vect{\gamma}(t),t) = \partial_t\scalar + \frac{\d \vect{\gamma}}{\d t}\cdot\nabla\scalar = 0\,.
\end{equation}
This means that along these points the solution stays constant in time and can be mapped back to its initial values:
\begin{equation}
    \scalar(\vect{\gamma}(t),t)=\scalar_0(\vect{x})\,.
\end{equation}
The curves $\gamma\colon\mathbb{R}_+\to\domain$ of constant values $\scalar$ are termed \define{characteristic curves}.
Evolving the solution of \cref{eq-def:characteristic_curve} for all $\vect{x}\in\domain$ simultaneously is done using the \define{characteristic map} $\vect{X}\colon\domain \times \mathbb{R}_+\to \domain$. Thus the characteristic map $\map(\vect{\gamma}(t),t)=\vect{\gamma}(0)=\vect{x}$ evolves all characteristic curves consecutively and the solution of \cref{eq-def:CMM-advect} at  time $t$ can be computed by the composition with the initial condition
\begin{equation}
    \scalar(\vect{x},t)=\scalar_0(\map(\vect{x},t))\,.
\end{equation}
Since $\map(\vect{x},t)$ traces back the characteristic curve to its initial position, it is called \define{backward characteristic map}. This formulation requires an assumption on the regularity of the velocity field as well as on its divergence (to prevent crossings of characteristic curves) so that we maintain a smooth one-to-one relabeling symmetry between Eulerian positions and Lagrangian particles, i.e. that the characteristic map is a diffeomorphism \citep{arnol2013mathematical}.
This is the case if the velocity field is smooth and divergence-free, i.e.,
\begin{equation}
    \nabla\cdot{\u} = 0\,.
\end{equation}

From the definition of the map $\map(\vect{\gamma}(t),t)=\vect{\gamma}(0)=\vect{x}$ it directly follows that $\frac{\d}{\d t} \map(\gamma(t),t) = 0$ and thus $\map$ solves
\begin{align}
\label{eq-def:CMM-map}
    \begin{cases}
     \partial_t\map + {\u}\cdot \nabla\map =0 & \text{for } (\vect{x},t)\in \domain\times\mathbb{R}_+\\
   \map(\vect{x},0) =\vect{x}           &  \vect{x}\in \domain\,.
    \end{cases}
\end{align}

Note, that \cref{eq-def:CMM-advect} is non-linear if ${\u}$ is dependent on the state $\scalar(\vect{x},t)$. Therefore, the characteristic map is dependent on the state and thus varies for different initial conditions.

\subsection{Advecting the Particle Distribution Function}
\label{subsec:advect-distr-func}

\begin{figure}[htp!]
    \centering
    \includegraphics[width=0.5\textwidth]{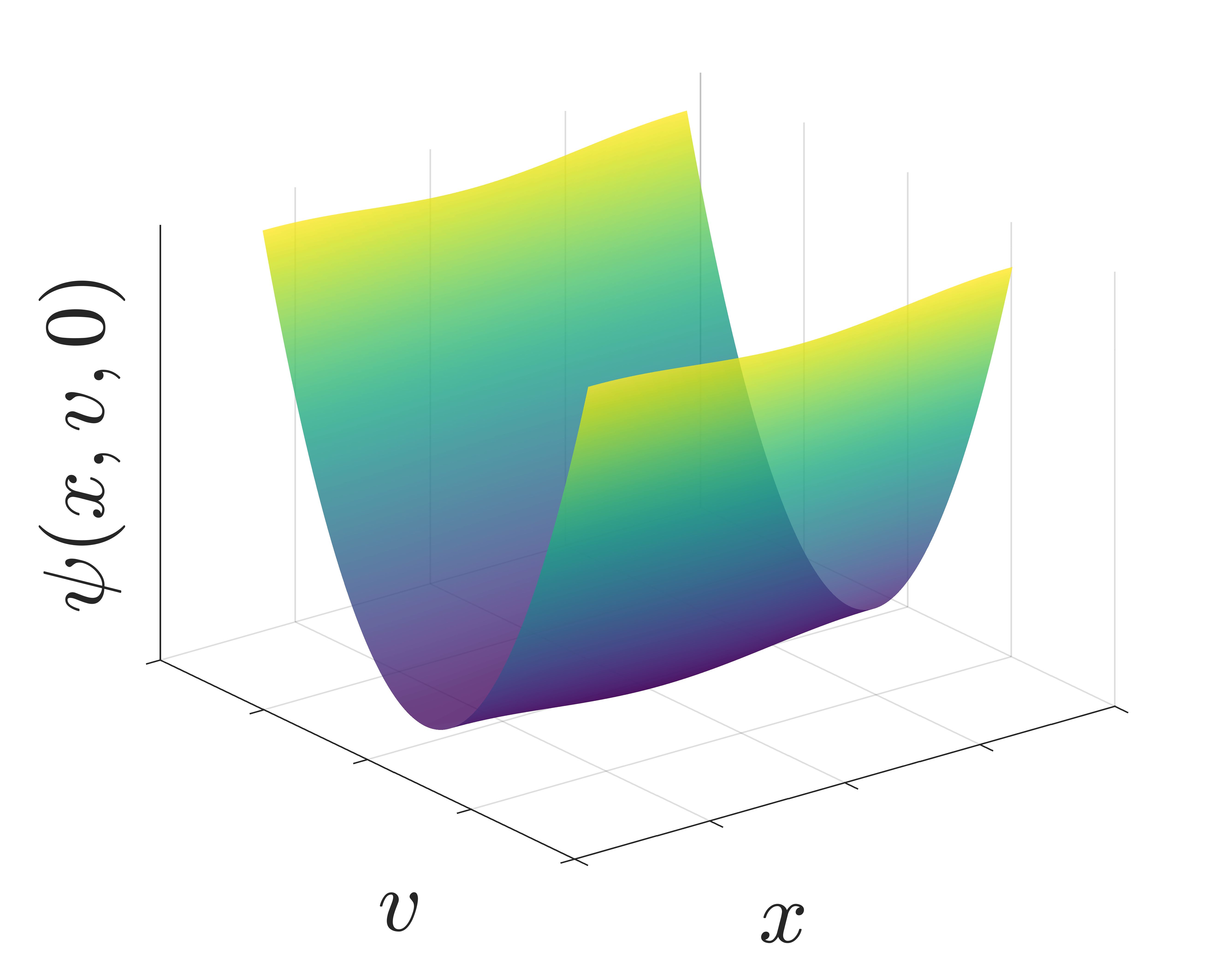}
    \caption{Stream function $\psi$ for the initial condition of the two-stream instability.}
    \label{fig:psi_2stream}
\end{figure}

Starting at \cref{eq:vlasov-advect} we can use the backward characteristic map $\vmap(x,v,t)=(X(x,v,t),V(x,v,t))$ to relate the distribution function to its initial condition:
\begin{equation}
    f(x,v,t) = f_0(X(x,v,t),V(x,v,t))\,,
\end{equation}
by solving
\begin{align}
    0 &=\partial_t \vmap + {\u}\cdot \nabla^\u \vmap, \quad\text{where}\quad {\u} = (v, \partial_x \phi)\\
    \partial_{xx}\phi &= l\int f\, \d v -1\,.
\end{align}
Since $\nabla^\u\cdot \vect{u} = 0$, we can define an associated 2D stream function
\begin{align}
    \psi(x,v,t) = \frac{v^2}{2}-\phi(x,t)\,,
\end{align}
such that ${\u} = \nabla^\u\times {\psi}\vect{e}_z$. For convenience $\psi$ is shown in \cref{fig:psi_2stream} for the two-stream instability at the initial condition. 
For implementation reasons of the CMM method, we assume periodic boundary conditions in $v$ additional to $x$.

\subsubsection{Periodic Boundary Conditions in Velocity Space}
\label{subsec:periodify}

\begin{figure}[htp!]
    \centering
    \setlength\figureheight{0.2\linewidth}%
    \setlength\figurewidth{0.9\linewidth}%
    \includetikz{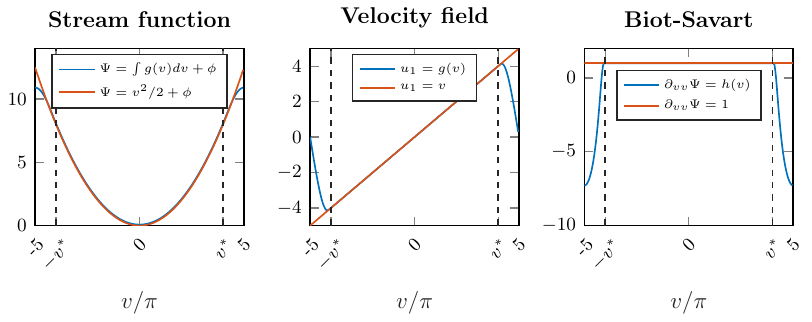}
    \caption{Periodization of the velocity field using a rescaled and shifted smooth Heaviside function $h(v)$. 
    From left to right: resulting periodized stream function, 
    a smooth continuation of the velocity component $u_1$ at different grid locations $v$, and the rescaled and shifted Heaviside function. }
    \label{fig:vel_periodify}
\end{figure}

The map advection assumes periodic boundaries in the $x$ and $v$ direction. 
For regularity of the map, we therefore have to maintain a periodic continuation of the velocity field $g(v)\colon \mathbb{R}\to\mathbb{R}$,
which is $g(v) = v$ for  $v \in \domain_{v*}:=[-v^*,v^*]\subset\domain_v$ and a smooth continuation from $-v^*$ to $v^*$ for $v\in\domain_v\setminus\domain_{v^*}$.
An example of such a function $g$ is visualized in the center of \cref{fig:vel_periodify}. 

We define $g(v)$ as the anti-derivative of $h(v):=\partial_v g(v)$.
As seen in the right of \cref{fig:vel_periodify} $h$ is a smooth heaviside function with $h(v)=1$ at $v\in\domain_{v^*}$.
For our purpose we choose: 
\begin{equation}
    \label{eq:heaviside}
    \sigma(v) = 1 - \frac{L_v}{a^2} \eta \left( \frac{\abs{v}-L_v}{a} \right) ,\quad \text{with}\quad a = 0.1L_v, 
\end{equation}
where 
\begin{equation}
    \eta (v)={\begin{cases} \frac{1}{I} \, e^{-1/(1-|v|^{2})}&{\text{ if }}|v|<1\\0&{\text{ if }}|v|\geq 1\end{cases}}\, ,
\end{equation}
where $I = \int e^{-1/(1-|x|^{2})}\d{x}$. Noting the fact that for a periodic domain the total circulation 
$\int_{\domain_x\times\domain_v}\vect{\nabla} \times \vect{u} \,\d x \d v$ of the associated velocity field $\vect{u}=(u_1,u_2) = (g(v),\partial_x \phi)$ vanishes. 
Thus from
 \begin{equation}
     \vect{\nabla}\times \vect{u} = \partial_v g(v) - \partial_{xx}\phi(x)
\end{equation}
we see that the total circulation can only vanish if $h(v)$ fulfills
\begin{equation}
    \label{eq:h_circ_condition}
    \int_{\domain_v} h(v) \d v=0\,.
\end{equation}
To ensure numerically that \cref{eq:h_circ_condition} and $h(v)=1$ in $v\in\domain_{v^*}$ we rescale and shift  $\sigma$ defined in \cref{eq:heaviside} such that:
\begin{equation}
    \label{eq:rescale_shifted}
    h(v) = \frac{\sigma - \int_{\domain_v} \sigma(v) \d v}{\max(\sigma - \int_{\domain_v} \sigma(v) \d v)}\,.
\end{equation}
The function $h$ is visualized together with the resulting constant velocity component $u_1$ and the associated stream function 
$\psi(x,v)$ along the $v$-direction in \cref{fig:vel_periodify}. It can be seen that the plotted functions $\Psi,g,h$ intersect with the exact values in the domain $\domain_{v^*}$.
The boundary layer in which the functions differ from the exact values is of size $a=L_v-v^*$. It is a tuning parameter of our algorithm, which will be discussed in the \cref{appdx:BLayer}.

This boundary issue is an inherent problem for mapping-based methods for kinetic equations. We denote by $\xmap_{[0, t]}$ the exact flow map of the unmodified phase space velocity $\vect{u}^e = (v, \partial_x \phi)$. We note that since $u_2 = \partial_x \phi$ is independent of $x$ and is in general non-zero, we have that for any strip $R = S^1 \times [v_0, v_1]$ the image under forward flow will generally exit the strip, i.e. $ \xmap_{[0, t]} (R) \not\subset R$ preventing the direct use of a bounded computational domain. However, the relevant error metric is the advection error of the density $f$, i.e. $\| f_0 (\vmap_{[t, 0]}) - f_0 (\xmap_{[t, 0]}) \|_{\infty} \leq C \| \nabla f_0 |_{\xmap_{[t,0]}} (\vmap_{[t, 0]} -\xmap_{[t, 0]}) \|_{\infty}$. Therefore, it is sufficient to control the backward map error where $\nabla f_0$ is nonzero. More generally, if one is allowed to assume a compactly supported initial condition, we let $R_0 = \text{supp}(f_0)$ the support of $f_0$ and $R_t$ the support of $f = f_0 \circ \xmap_{[t,0]}$ at time $t$. Writing $\mathbf{1}_{R}$ for the indicator function of the set $R$, we have that $\mathbf{1}_{R_t} = \mathbf{1}_{R_0} \circ \xmap_{[t, 0]}$. The advection error then becomes 
\begin{align}
\| f_0 (\vmap_{[t, 0]}) - f_0 (\xmap_{[t, 0]}) \|_{\infty} & \leq C \| \nabla f_0 \|_{\infty} \| (\vmap_{[t, 0]} -\xmap_{[t, 0]}) \mathbf{1}_{R_0} \circ \xmap_{[t,0]} \|_{\infty} \\
 & = C \| \nabla f_0 \|_{\infty} \| (\vmap_{[t, 0]} -\xmap_{[t, 0]}) \mathbf{1}_{R_t} \|_{\infty} .
\end{align}

If we further suppose that there exists some strip $\Omega_{v^*} = S^1 \times [-v^*, v^*]$ such that $R_t \subset \Omega_{v^*}$ for all $t \in [0, T]$, we have
\begin{equation}
    0 = (\partial_t + \vect{u} \cdot \nabla ) \vmap_{[t,0]} - (\partial_t + \vect{u}^e \cdot \nabla ) \xmap_{[t,0]} = (\partial_t + \vect{u}^e \cdot \nabla ) (\vmap_{[t,0]} - \xmap_{[t,0]}) 
\end{equation}
inside $\Omega_{v^*}$ so that $(\vmap_{[t,0]} - \xmap_{[t,0]})$ satisfies the advection equation under $\vect{u}^e$ in $\Omega_{v^*} \times [0, T]$ with initial condition $\vmap_{[0,0]} - \xmap_{[0,0]} = 0$ in $\Omega_{v^*}$. It follows from the advection property that this difference remains zero in the set $\{ (x, v, t) \in \Omega_{v^*} \times [0, T] | \xmap_{[t, 0]}(x, v) \in \Omega_{v^*} \} \supset \text{supp}(f)$. Hence $\| f_0 (\vmap_{[t, 0]}) - f_0 (\xmap_{[t, 0]}) \|_{\infty} = 0$.

To apply this flow modification, we must know \emph{a priori} that $\text{supp}(f) \subset \Omega_{v^*}$ for all $t \in [0, T]$. If the support eventually exits the strip, a reinitialization is then required.

\newcommand{\CoarseGrid}{G}
\newcommand{\Xf}{G_\mathrm{f}}
\newcommand{\XPsi}{G_\psi}
\subsection{Numerical Implementation}

The implementation details closely follow \cite{yin2021characteristic} and are based on the gradient-augmented level set method (GALS) and Jet-scheme frameworks. Here we only review the most relevant parts. However, for a detailed explanation, we refer the reader to \cite{yin2021characteristic,yin2021thesis}.

\paragraph{Spatial Discretization.}

The CMM implementation makes use of three different kinds of grids, which are summarized in \cref{tabl:grids}.

\begin{table}[htp!] 
\centering                
\begin{tabular}{lll}
\toprule
Grid                  & Size       & Explanation\\
\midrule
Coarse grid $\CoarseGrid$   & $\Ncoarse\times \Ncoarse$ &Stores solution of advection equation \\
Sampling grid $\Xf$     & $\Nfine \times \Nfine$       &Stores the samples of the initial condition\\
Velocity grid   $\XPsi$    & $\Npsi \times \Npsi$          &Velocity samples\\
\bottomrule
\end{tabular}
\caption{Overview of the different computational grids used in this work.}
\label{tabl:grids}
\end{table}

The computational domain is discretized with an equidistant grid $\CoarseGrid$. 
In the initial step of the simulation, the PDF is sampled on the sampling grid $\Xf$ of size $\Nfine\times \Nfine$. Therefore the map is interpolated using cubic Hermite interpolation and the PDF is evaluated $f=f_0\circ\mathcal{H}_{\Xf}[\vmap]$ at the corresponding positions.
Next, the electric potential is computed on the sampling grid, by integrating over $f$ in space and solving the 1D Poisson equation:
\begin{equation}
    \label{eq:poisson}
    \widehat\phi(k_x,t) =\frac{1}{k_x^2} \, \mathcal{F}[1-\int f(\cdot,v,t) \d v](k_x)\,,
\end{equation}
where $\mathcal{F}$ denotes the Fourier transform (FFT) defined in \cref{sec:pseudo-spectral}.
Depending on the resolution of the velocity grid $\XPsi$ we perform zero padding, to upsample the solution from the sampling grid $\Xf$ to the velocity grid $\XPsi$.
From the solution of \cref{eq:poisson}, we form the stream function:
\begin{equation}
    \psi(x,t) = g(v) - \mathcal{F}^{-1}[\widehat\phi(\cdot,t)](x)
\end{equation}
with the precomputed periodic velocity field $g(v)$.
The Hermite representation $\mathcal{H}_{\XPsi}[\psi]$ of the stream function $\psi$ involves the derivatives of $\psi$ with respect to $x$ and $v$ direction, as well as the cross derivative $xv$. Here we make use of $h(v)=\partial_v g(v)$ and $\mathcal{F}^{-1}[\hat\phi(k_x,t)](x)$ that have been already computed. Thus, no further 2D FFT is necessary since the cross derivative vanishes.
From the obtained stream function in Hermite form we can compute the velocity using the identity:
\begin{align}
    \vect{u}=\nabla \times \mathcal{H}_{\XPsi}[\psi].
\end{align}

\paragraph{Temporal Discretization.}

With the given velocity field we integrate \cref{eq-def:characteristic_curve} using a 3rd order Runge--Kutta backwards in time.

\paragraph{Compositional Structure and Remapping Criteria.}

The incompressibility of the velocity field implies that the map is volume preserving, which means that the determinant of its Jacobian equals one. Numerically this property is violated and remapping becomes necessary. The compositional structure of the map, i.e. its decomposition into sub-maps, is hereby helpful and increasingly complicated spatial features in the solution can thus be adaptively resolved. To this end, we monitor the determinant of the Jacobian of the map over time, which should be equal or close to one. We introduce the following error criterion \rev{to approximate the} a posteriori error estimate \citep{yin2021characteristic}:
\begin{equation}
     e_\text{det}^{n} \, = \, || \det \nabla \vmap^n  - 1 ||_{\infty} \, \le \, \delta_\text{det}\,,
     \label{eq:remapping}
\end{equation}
where $\delta_\text{det}$ is the threshold for the Jacobian error, also called \define{imcompressibility threshold} in the forthcoming. 
In case that $e_\text{det}^{n}$ exceeds this threshold, we store the sub-map $\vmap^n$\rev{, reinitialize} the new map with the identity and proceed with the time evolution.
This remapping procedure guarantees volume preservation of the map up to a certain accuracy.
Therefore $\delta_\text{det}$ is chosen as small as possible, but large enough, such that all sub-maps can be stored in the active memory. For our fine-scale simulations, we discuss the influence of $\delta_\text{det}$ in \cref{appx:fine-scale-stats}.

\subsection{Error Analysis}
\label{subsec:ErrorAnalysis}

The numerical method for the evolution of the characteristic map is the same as the one studied in \cite{yin2021characteristic} with a key difference in the definition of the velocity field. In \cite{yin2021characteristic}, the numerical velocity field $\tilde{\u}$ is obtained from the map $\vmap$ using the Biot-Savart law on the vorticity field: $\tilde{\u} = - \nabla \times \Delta^{-1} \omega_0 (\vmap ( {\bm x}, t))$. In the case of Vlasov--Poisson, the $x$-component of the velocity field is given analytically and is exact within the region $[-v^*, v^*]$. The $v$-component of the velocity field is obtained from solving a 1D Poisson equation for the advected distribution function:
\begin{equation}
    u_v(x) = \partial_x \phi = l \partial_x \Delta_x^{-1} \left(  \int f_0 (\vmap(x, v, t) ) \d v -1  \right).
\end{equation}

The error term for the velocity field is given by
\begin{equation}
    u_v - u^*_v \approx l \partial_x \Delta_x^{-1} \left(  \int \nabla f_0 \cdot (\vmap(x, v, t) - {\vmap}^*(x, v, t) ) \d v \right),
\end{equation}
for smooth initial condition $f_0$. We note that the velocity error is a $1^{st}$ anti-derivative of the map error in the $x$ component and hence can be bounded by the map error. This is enough to close the inequalities in the error bounds as seen in \cite{yin2021characteristic} and \cite{yin2023characteristic}. Therefore, we get the same error estimates as in \cite{yin2021characteristic}.

 Using Hermite cubic interpolation coupled with $3^{rd}$ order in time local Lagrange extrapolation and $3^{rd}$ order Runge--Kutta integrators for the one-step maps, we expect global $3^{rd}$ order for the map error:
 \begin{equation}
     \| \vmap(\cdot, t^n) - \vmap^n \|_{\infty} = \mathcal{O} (\Delta x^3 + \Delta t^3).
 \end{equation}
Assuming smooth initial condition, this also implies a $3^{rd}$ order convergence for the distribution function $f$. For mass conservation, we have
\begin{equation}
    \mathcal{M} = \iint f_0  \d x \d v = \iint f_0 \circ \vmap^n \det (\nabla \vmap^n) \d x \d v
\end{equation}
by a change of variables for the integration. It follows that the mass conservation error is bounded by the minimum of the distribution function error and the volume conservation error. The conservation errors for momentum and total energy are also controlled by the distribution function error and are therefore $3^{rd}$ order. 
In summary\rev{,} we define the following errors:
\begin{alignat}{3}
    \Delta \vmap &=  \Vert\vmap^\text{ref}(\cdot,t^n)- \vmap^n \Vert_\infty    &\quad& \text{(map error)}\,,                      \label{eq:map_error}\\
    \Delta f &=  \Vert f^\text{ref}(\cdot,t^n)- f^n \Vert_\infty               &\quad& \text{(distribution function error) }\,,    \label{eq:dist_error}\\
    \Delta \mass &=  \abs{\mass(t^n)- \mass(t^0)}                   &\quad& \text{(mass conservation error)} \,,       \label{eq:MconsError}\\
    \Delta \momentum &=  \abs{ \momentum(t^n)- \momentum(t^0) }     &\quad& \text{(momentum conservation error)}\,,    \label{eq:PconsError}\\
    \Delta \Etot &=  \abs{ \Etot(t^n)- \Etot(t^0) }                 &\quad& \text{(energy conservation error)}\,,      \label{eq:EconsError}\\
\end{alignat}
that are of $3^{rd}$ order and will be studied numerically in \cref{sec:numtest}.

\section{Pseudo-Spectral Discretization of Vlasov--Poisson}
\label{sec:pseudo-spectral}

Assuming the same periodized
velocity field as introduced in \cref{subsec:periodify} we can utilize a classical Fourier pseudo-spectral method \cite{HussainiYousuffZang1987}.
Such a periodization approach where a known profile has been substrated and the perturbation becomes periodic has been used in \cite{FroehlichSchneider1993} in the context of combustion problems.
The Fourier projection of the density distribution function $f$ is given by
\begin{equation}
P_N f (\bm x) = f_N (\bm x) = \sum_{|{\bm k}| \lesssim N/2} \widehat f_k \, e^{i {\bm k} \cdot {\bm x}} \; , \; \widehat f_{\bm k} = \frac{1}{(2 \pi)^d} \int_{\mathbb{T}^d}\, f({\bm x}) \, e^{-i {\bm k}  \cdot {\bm x}} \, d{\bm x}\,
\label{eq:Fourierprojector}
\end{equation}
where $ \mathbb{T} = \mathbb{R} / (2 \pi \mathbb{Z})$ and $d$ \rev{is} the dimension.
Note that $|{k}| \lesssim N/2$ is understood in the sense $-N/2 \le k < N/2$ and correspondingly in higher dimensions for each component of $\bm k$.

Applying the spectral discretization to the one+one-dimensional Vlasov--Poisson equation ($d=2$) 
with periodic boundary conditions and suitable initial condition $f({\bm x},t=0) = f_0({\bm x})$ yields the Galerkin scheme
\begin{equation}
\partial_t f_N + \left( P_N ({\u_N}\cdot \nabla^{\u} f_N) \right) \, = \, 0  \quad {\text{for}} \quad {\bm x} \in  \mathbb{T}^2 \quad {\text{and }} \quad t>0 ,
\end{equation}
which corresponds to a nonlinear system of $N\times N$ coupled ODEs for $\widehat f_{\bm k}(t)$ with $|{\bm k}| \lesssim N/2$.
A pseudo-spectral evaluation of the nonlinear term is utilized, and the product in physical space is fully dealiased. In other words, the Fourier modes retained in the expansion of the solution are such that $|{\bm k}| \le k_C$, where $k_C$ is the desired cut-off wave number, but the grid has $N= 3k_C$ points in each direction, versus $N= 2k_C$ for a non-dealiased, critically sampled product. This dealiasing makes the pseudo-spectral scheme equivalent to a Fourier--Galerkin scheme up to round-off errors \cite{HussainiYousuffZang1987}, and is thus conservative for the $L_2$ norm. 

For time discretization of the resulting ODE systems we use a classical Runge--Kutta scheme of order 3.
For details on the convergence and stability of the above spectral schemes, we refer to \cite{bardos2015stability}.

\medskip

Note that \cref{eq:Vlasov--Poisson} can be efficiently solved in Fourier space assuming periodic boundary conditions.
Therefore we discretize \cref{eq:vlasov-advect} in space using pseudo-spectral scheme \cite{HussainiYousuffZang1987}
\begin{align}
    \hat\phi(k_x,t) &= \mathcal{F}[\phi(\cdot ,t)](k_x)=\int \phi(x,t)\ e^{-i k_x x}\,\d x,\\
    \hat f(k_x,k_v,t) &=\mathcal{F}[f(\cdot,\cdot ,t)](k_x,k_y)=\int \int f(x,v,t)\ e^{-i (k_x x+k_v v)}\,\d x \d v,\\
\end{align}
and the inverse:
\begin{align}
    \phi(x,t) &=\frac{1}{2\pi}\int\hat \phi(k_x,t)\ e^{i k_x x}\,\d k_x,\\
    f(x,v,t) &=\frac{1}{4\pi^2}\int \int \hat f(k_x,k_v,t)\ e^{i (k_x x+k_v v)}\,\d k_x \d k_v.\\
\end{align}
Setting this into \cref{eq:Vlasov--Poisson} yields
\begin{align}
    \hat\phi(k_x,t) &= \frac{1}{k^2_x} \mathcal{F}[(1 - l \int f\, \d{v})](k_x)=\frac{1}{k^2_x}\left(\mathcal{F}[1](k_x) - l  \underbrace{\mathcal{F}[\int f\, \d{v})](k_x)}_{\hat f(k_x,0,t)}\right)\\
    &=\frac{1}{k_x^2}(\delta(k_x) - l\hat f(k_x,0,t)) \,.  \label{eq:Vlasov--Poisson-fft}
\end{align}
After solving \cref{eq:Vlasov--Poisson-fft} we solve \cref{eq:vlasov-advect} \rev{with a} pseudo-spectral \rev{method}. To avoid aliasing errors we remove 1/3 of the Fourier modes, i.e., the high frequencies. This filtering is known as the 2/3 rule and has proven to lead to the conservation of the $L_2$ norm of the solution \cite{bardos2015stability}. Furthermore, we remark that the total mass $\mass$ is conserved over time since it is equal to the first Fourier coefficient that is unchanged during time integration.

\section{Numerical Studies} 
\label{sec:numtest}
In the following, we study Landau damping in the linear and non-linear regime and the two-stream instability. Convergence tests are carried out and the CMM is assessed and benchmarked with respect to the pseudo-spectral method and results in the literature. 

With this publication, we provide an open-source MATLAB code as well as a C\texttt{++} CUDA single GPU implementation \footnote{Code repository \url{https://github.com/orgs/CharacteristicMappingMethod/repositories}}.
The MATLAB repository also includes the implementation of the pseudo-spectral code that we use for the convergence analysis. The MATLAB convergence tests are computed on a \texttt{Intel 11th Gen i7-11850H} 8-core CPU platform.

In addition to the convergence tests, we run \rev{high resolution} simulations for all three test cases to access \rev{the fine-scale features of the solution}. The simulations are performed with a CUDA implementation that is run on \texttt{NVIDIA A100 SXM4} 80 GB GPU. 

We define experimental order of convergence (EOC) for two consecutive simulations with coarse map $N_c=2N$ and $N_c=N$ as:
\begin{equation}
    \mathrm {EOC}(N)={\frac {\log {\frac {\|f-f_{2N}\|}{\|f-f_{N}\|}}}{\log {2}}} 
    \label{eq-def:EOC}
\end{equation}
for a given norm $\Vert{\cdot}\Vert$.
\subsection{Landau Damping}
In this section we study Landau damping initialized with a perturbed Maxwellian distribution:
\begin{equation}
\label{eq:landau_damping_inicond}
    f_0(x,v) =  \left(1+\epsilon\cos(kx) \right)\frac{1}{\sqrt{2\pi}} \, e^{-v^2/2}\,.
\end{equation}
All relevant physical parameters are stated in \cref{tabl:landaudamping}.
Two perturbation regimes are studied: 
\begin{enumerate}
\item linear Landau damping regime with $\epsilon=0.05$
\item non-linear Landau damping regime with $\epsilon=0.5$
\end{enumerate}

\begin{table}[htp!]
    \centering
    \begin{tabular}{l l}
    \toprule
       name  & value \\
       \midrule
         \rev{spatial domain}    &     \rev{$\domain_x\in[0, 2\pi/k]$}\\
         \rev{velocity domain} & \rev{$\domain_v\in[-4\pi,4\pi]$}\\
         \rev{effective vel. domain} &\rev{$\domain_{v^*}\in[-3.8\pi,3.8\pi]$}\\
         wave number $k$ &    $0.5$\\
         $\epsilon$ (linear)     &   $5 \times 10^{-2}$ \\ 
         $\epsilon$ (non linear)     &   $5 \times 10^{-1}$ \\ 
         \bottomrule
    \end{tabular}
    \caption{Parameters taken from \cite{PhamHelluyCrestetto2013} for the Landau damping.}
    \label{tabl:landaudamping}
\end{table}

\subsubsection{Linear Landau Damping}
\label{subsec:linear_landau_damping}
We initiate our analysis by conducting a convergence assessment of our CMM algorithm, employing the simulation parameters outlined in \cref{tabl:landaudamping}.

\paragraph{Convergence Tests.}
To evaluate the spatial and temporal convergence properties of our proposed numerical algorithm, we vary the discretization parameters detailed in \cref{tabl:convergence-params}.
\begin{table}[htp!]
\centering
 \caption{Grid parameters of the convergence tests.}
 \begin{subtable}{.48\linewidth}
    \begin{tabular}{l l}
    \toprule
       name  & value \\
       \midrule
         $\Ncoarse$ coarse grid $\CoarseGrid$   &   32, 64, 128, 256, 512\\
         $\Nfine$ sampling grid & 1024\\
         $\Npsi$ velocity grid & 1024\\
         $\Delta t$ time step size &   L/1024\\
         $\Ncoarse^\text{ref}$ ref grid & 1024\\
         \bottomrule
    \end{tabular}%
    \caption{Convergence tests in space}%
    \label{tabl:convergence_space}%
    \end{subtable}
    \hfill
     \begin{subtable}{.48\linewidth}
    \begin{tabular}{l l}
    \toprule
       name  & value \\
       \midrule
         $\Ncoarse$ grid resolutions $\CoarseGrid$   &  1024\\
         $\Nfine$ sampling grid & 1024\\
         $\Npsi$ velocity grid & 1024\\
         $\Delta t$ time step sizes  &   $2^{-7},2^{-6},2^{-5},2^{-4}$\\
         ref time step size & $2^{-8}$\\
         \bottomrule
    \end{tabular}
    \caption{Convergence tests in time}
    \label{tabl:convergence_time}
    \end{subtable}
    \label{tabl:convergence-params}
\end{table}

For the spatial convergence evaluations, we conduct simulations with varying spatial resolutions for the characteristic map $\vmap$, integrating them in time up to $T=10$, while keeping all other parameters constant. It is worth noting that we disable remapping during these convergence tests.

Upon reaching the final time point of $T=10$, we calculate errors as specified in \cref{eq:dist_error,eq:map_error,eq:MconsError,eq:PconsError,eq:EconsError} on the fine grid $\Xf$. As a reference, we use the simulation performed with $\Ncoarse=1024$. This self-comparison is done to assess the error associated with the map $\vmap$.

\begin{figure}[htp!]
    \centering
    \begin{subfigure}{0.46\linewidth}
    \setlength\figureheight{0.8\linewidth}%
    \setlength\figurewidth{0.8\linewidth}%
    \caption{}
    \includetikz{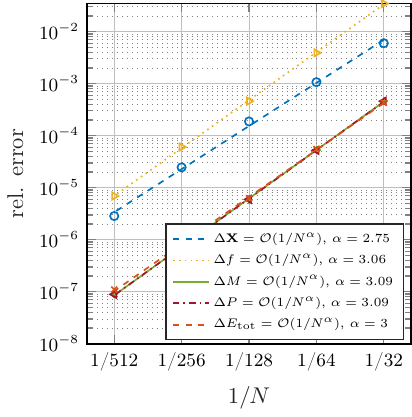}
    \label{fig:landau_convergence_space}
    \end{subfigure}
    \hfill
    \begin{subfigure}{0.46\linewidth}
    \caption{}
    \setlength\figureheight{0.8\linewidth}%
    \setlength\figurewidth{0.8\linewidth}%
    \includetikz{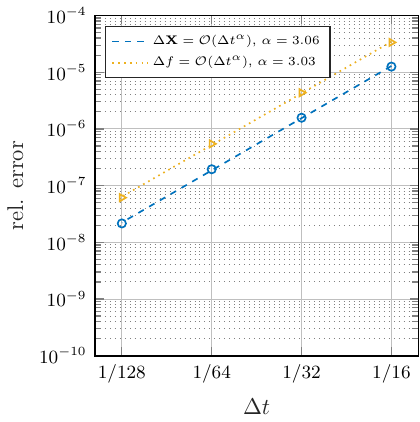}
    \label{fig:landau_convergence_time}
    \end{subfigure}
    \caption{Convergence studies for the linear Landau damping test case at time $t=10$ with physical parameters listed in \cref{tabl:landaudamping}. Shown are the relative errors that correspond to the error measures defined in \cref{eq:dist_error,eq:map_error,eq:MconsError,eq:PconsError,eq:EconsError}.\\
    \textbf{(a)} Spatial convergence of evolutionary quantities $f$, $\vmap$ and the conserved quantities $\Delta \mass,\Delta \momentum \Delta \Etot$ using the parameters listed in \cref{tabl:convergence_space}.\\
     \textbf{(b)} Temporal convergence with parameters stated in \cref{tabl:convergence_time}.
    }
    \label{fig:landau_convergence}
\end{figure}

As elaborated in \cref{subsec:ErrorAnalysis}, we anticipate third-order convergence in space for both the $L_\infty$ error associated with the map and the Probability Density Function (PDF), as well as third-order conservation errors. These expectations are substantiated by our numerical experiments, as illustrated in \cref{fig:landau_convergence}, and the experimental convergence order provided in \cref{table:EOC-landaudamping}.

\begin{figure}[htp!]
    \centering
    \setlength\figureheight{0.3\linewidth}%
    \setlength\figurewidth{0.4\linewidth}%
    \includetikz{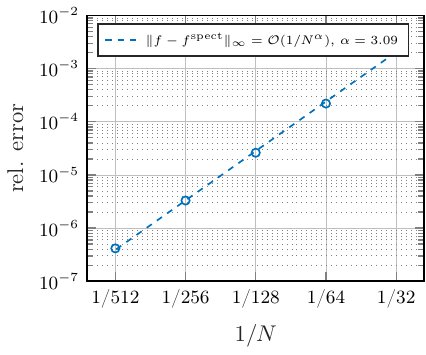}
    \caption{Spatial convergence of the linear Landau damping test case. Relative error as a function of $N$ of CMM with respect to the Fourier--Galerkin reference solution computed at resolution $512^2$.}
    \label{fig:landau_convergence_fourier}
\end{figure}

Furthermore, we can affirm third-order spatial convergence when comparing the proposed method to a reference solution of a pseudo-spectral computation conducted on a computational grid with dimensions of $512\times512$. The results are succinctly summarized in the final column of \cref{table:EOC-landaudamping} and visually represented in \cref{fig:landau_convergence_fourier}.
\begin{table}[htp!]                                                                                  
\centering                                                                                     
\begin{tabular}{ccccccc}                                                                        
\toprule                                                                                       
$1/N$ & $\Delta f$ & $\Delta \mathbf{X}$ & $\Delta M$ & $\Delta P$ & $\Delta \Etot$ & $\norm{f-f^\text{spectr}}_\infty$\\
\midrule                                                                                       
1/32 & - & - & - & - & - \\                                                                    
1/64 & 3.1 & 2.5 & 3.2 & 3.1 & 3.0  &3.4\\                                                          
1/128 & 3.1 & 2.5 & 3.1 & 3.1 & 3.1 &3.1\\                                                         
1/256 & 2.9 & 2.9 & 3.1 & 3.1 & 3.0 &3.0\\                                                         
1/512 & 3.1 & 3.1 & 3.0 & 3.0 & 2.8 &3.0\\                                                         
\bottomrule                                                                                    
\end{tabular}                                                                                  
\caption{Spatial experimental order of convergence EOC as defined in \cref{eq-def:EOC} for the linear Landau damping test case.}                                                                       
\label{table:EOC-landaudamping}                                                                     
\end{table}   

Similarly, we observe third-order temporal convergence, as evidenced by the results showcased in \cref{fig:landau_convergence_time} and detailed in \cref{table:EOC-time}. This assessment of temporal convergence is achieved by comparing simulations conducted with various time steps integrated up to $T=10$. Importantly, for the temporal convergence examinations, we maintain a consistent spatial resolution as outlined in \cref{tabl:convergence_time}.

\begin{table}[htp!]                                                                                       
\centering                                                                                          
\begin{tabular}{ccc}                                                                             
\toprule                                                                                            
$\Delta t$ & $\Delta f$ & $\Delta \mathbf{X}$  \\
\midrule                                                                                            
1/16 & - & - \\                                                                         
1/32 & 2.9 & 3.0  \\                                                             
1/64 & 3.0 & 3.0  \\                                                             
1/128 & 3.2 & 3.2 \\                                                           
\bottomrule                                                                                         
\end{tabular}                                                                                       
\caption{Temporal order of convergence EOC as defined in \cref{eq-def:EOC} for the Landau damping test case.}                                                                            
\label{table:EOC-time}                                                                          
\end{table}

\paragraph{Damping Rate, Sub-grid Resolution and Recurrence Time.}
A major feature of our method is that the local time stepping is not restricted by the spatial representation of the solution itself. Therefore fine-scale spatial properties can be transported with a larger time step to high accuracy without losing stability. This fact can be exploited to accurately predict the damping rate of linear Landau damping even for highly frequent initial perturbations \cref{eq:landau_damping_inicond}. 
One challenge that occurs in regular grid-based, periodic-in-$x$ numerical schemes of Vlasov Poisson is to overcome the \textit{velocity space filamentation problem}.
This is a well-known problem in numerical plasma simulations reported since the 1980s \cite{Klimas1987,ChenKnorr1976}. 
Here, the presence of minor initial spatial perturbations gives rise to oscillations in the velocity distribution, and the frequency of these oscillations increases as time progresses. Consequently, at a certain stage, the velocity grid becomes insufficient to accurately capture these increasingly rapid oscillations. This leads to an aliasing effect, wherein the numerical computation of charge density interprets high-frequency oscillations as if they were of lower frequency. The time at which this artifact comes to play was first discovered 1976 in \cite{ChenKnorr1976} and is known as the \textit{recurrence time} $\trec=2\pi/(k\Delta v)$.
Mehrenberger et al. showed in \cite{MehrenbergerNavoretPham2020} that $\trec$ is linked to the accuracy of the velocity quadrature when computing the charge density in the Poisson equation \cref{eq:Vlasov--Poisson}. 
This suggests integrating \cref{eq:Vlasov--Poisson} on a finely resolved grid, which however limits the time step size for conventional numerical schemes.

\begin{figure}[htp!]
    \centering
    \begin{subfigure}{0.6\linewidth}
    \centering
    \setlength\figureheight{0.7\linewidth}%
    \setlength\figurewidth{0.7\linewidth}%
    \caption{}%
    \includetikz{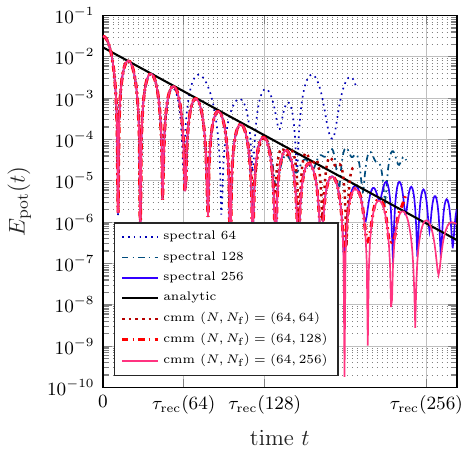}
    \label{subfig:ldamping_Epot_compare}
    \end{subfigure}\\
     \begin{subfigure}{0.6\linewidth}
     \centering
    \setlength\figureheight{0.7\linewidth}%
    \setlength\figurewidth{0.7\linewidth}%
       \caption{}%
       \includetikz{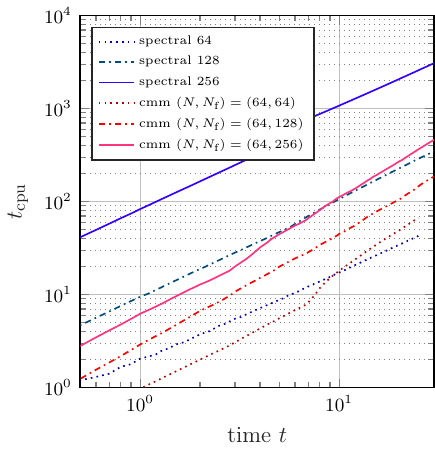}
       \label{subfig:ldamping_Epot_compare_cputime}
       \end{subfigure}
    \caption{Comparison of the damping of the potential energy for $k=0.5$.
    (a). Comparison of CPU-time for CMM and the spectral method for different grids (b).}
    \label{fig:landau_damping_epot_spectral_vs_CMM}
\end{figure}

As we show in \cref{fig:landau_damping_epot_spectral_vs_CMM}, we can overcome this problem when using a high-resolved sampling grid and a coarse map grid. 
\Cref{fig:landau_damping_epot_spectral_vs_CMM} compares the Landau damping for the potential energy at $k=0.5$ using pseudo-spectral and CMM simulations.
For both methods the recurrence time $\trec$ is clearly visible in \cref{subfig:ldamping_Epot_compare} at the time when the exponential damping of the potential energy saturates. With decreasing lattice spacing the recurrence time decreases as expected in the pseudo-spectral code. 
We note that increasing $\Ncoarse$ in CMM has no effect on the recurrence time.
However as expected, when the resolution of the sampling grid is increased, the recurrence time increases. Furthermore in general the recurrence time seems to be larger for the CMM code compared to the spectral code when using the same resolution. 

Note that the complexity for solving \cref{eq:Vlasov--Poisson} is the same in the pseudo-spectral and CMM case if the sampling grid and pseudo-spectral resolution are identical. 
However, since the time steps are not restricted by the resolution of the sampling grid, the CMM method saves several orders of magnitude when comparing the CPU time to our pseudo--spectral code, shown in \cref{subfig:ldamping_Epot_compare_cputime}. This is especially visible for increasing $N=\Npsi$. Therefore, the savings are purely contributed to solving the advection part.

\begin{figure}[htp!]
    \centering
    \setlength\figureheight{0.35\linewidth}%
    \setlength\figurewidth{0.35\linewidth}%
    \includetikz{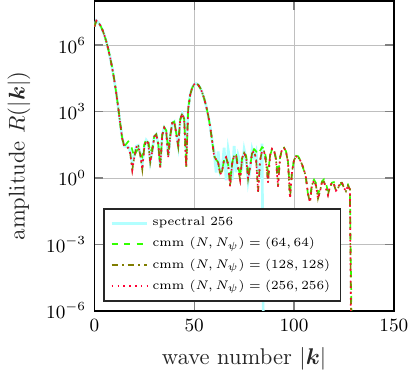}
    \caption{Radial Fourier spectra $R(\abs{\vect{k})}) = \int_{\abs{\vect{k}}-\Delta k}^{\abs{\vect{k}}+\Delta k}\abs{\mathcal{F}[f]({\vect{k})}}\d \vect{k}$ of the Fourier transformed snapshot $f(x,v,t)$ at time $t=25$ for the Landau damping test case. The coarsest CMM resolution can recover the finest spectral resolution.  }
    \label{fig:landau_damping_subgrid_spectral_vs_CMM}
\end{figure}

The CMM method provides an elegant way of circumventing the filamentation problem for numerical simulations, because of its sub-grid properties. As shown in the radial Fourier spectra in \cref{fig:landau_damping_subgrid_spectral_vs_CMM} the spectral content of the solution evaluated at a resolution of $256\times256$ at time $t=25$ can be resolved even with a coarse map grid of size $64 \times 64$. Note that the spectrum for the spectral simulation is cut at 85 because of the 2/3 rule.

\begin{figure}[hbp!]
    \centering
    \setlength\figureheight{0.3\linewidth}%
    \setlength\figurewidth{0.3\linewidth}%
    \includetikz{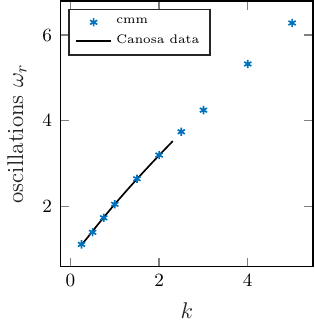}
    \includetikz{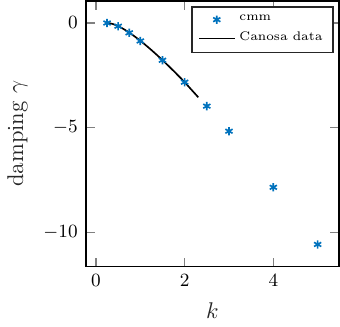}
    \caption{Comparison of the fitted damping rates and oscillation frequencies to the Canosa root finding data \cite{Canosa1973}. Here we use $ (\Ncoarse,\Npsi)=(64,1024)$, $\delta_\text{det}=0.05$.}
    \label{fig:landau_damping_rate_freq}
\end{figure}

As a result of the sub-grid property, we can resolve the finest oscillations in the solution, i.e. we can extract damping rates $\gamma$ and oscillation frequencies $\omega_r$ of the Landau damping even for large $k$. To obtain the damping rates we fit the peaks of the potential energy. We include all peaks after the first peak at $t=0$ that are above a certain threshold, here $\Epot(t)\ge 10^{-17}$. From the distance of the peaks, we calculate the mean frequency, which is shown in \cref{fig:landau_damping_rate_freq} together with the damping rate. We compare our results to the Canosa root finding data \cite{Canosa1973} as it is reported to yield more accurate results than other approximations \cite{finn2023numerical}.
In fact \cite{Canosa1973} claims an accuracy to the 5th digit. 
The quantitative results and relative errors compared to \cite{Canosa1973} are reported in \cref{table:damping_rate}. As indicated by the relative errors in \cref{table:damping_rate}, we obtain the most accurate results for large values of $k$. This error behavior was also observed in \cite{finn2023numerical}.

\begin{table}[htp!]                                                                       
\centering                                                                          
\begin{tabular}{ccccc}                                                              
\toprule                                                                            
$k$ & $\gamma$ & rel. err & $\omega_r$ & rel. err \\                                
\midrule                                                                            
0.2500 & -0.0032 & 0.4530 & 1.1062 & 0.0005 \\                                      
0.5000 & -0.1569 & 0.0233 & 1.3891 & 0.0187 \\                                      
0.7500 & -0.4726 & 0.0231 & 1.7283 & 0.0051 \\                                      
1.0000 & -0.8609 & 0.0113 & 2.0433 & 0.0013 \\                                      
1.5000 & -1.7812 & 0.0031 & 2.6400 & 0.0029 \\                                      
2.0000 & -2.8274 & 0.0001 & 3.1948 & 0.0018 \\                                      
\bottomrule                           
\end{tabular}                         
\caption{Numerical values of the Landau damping rate $\gamma$\rev{,} oscillation frequency $\omega_r$ and their relative error compared to \cite{Canosa1973} for $(\epsilon,k)=(0.05,5)$.}
\label{table:damping_rate} 
\end{table}  

\begin{figure}[hbp!]
    \centering
    \begin{subfigure}{0.48\linewidth}
        \centering
        \setlength\figureheight{0.7\linewidth}%
        \setlength\figurewidth{0.8\linewidth}%
        \caption{}
        \includetikz{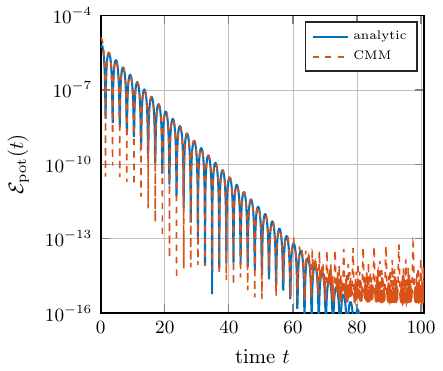}
        \label{fig:landau_damping_epot_non_lin}
    \end{subfigure}
    \begin{subfigure}{0.48\linewidth}
    \centering
        \setlength\figureheight{0.7\linewidth}%
        \setlength\figurewidth{0.8\linewidth}%
        \caption{}
        \includetikz{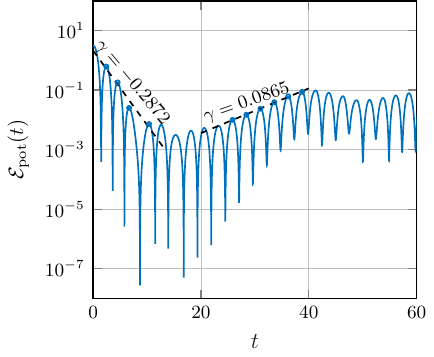}
        \label{fig:landau_damping_epot_non_lin}
    \end{subfigure}
    
    \caption{Potential energy of linear ($\epsilon = 0.001$) and non-linear ($\epsilon = 0.5$) Landau damping at $k=0.5$ using a fine grid computation $(\Ncoarse,\Nfine,\Npsi) = (1024, 4096, 4096)$. \\
    a) Linear Landau damping resolved down to machine precision without recurrence phenomena.\\
    b) Non-linear Landau damping and the different decay rates. The decay rates are determined by fitting the maxima that are plotted as circles.}
    \label{fig:Landau_damping_without_recurence}
\end{figure}

Lastly, we show in \cref{fig:Landau_damping_without_recurence} that we can resolve the potential energy until machine precision when the sampling grid is increased massively. This computation utilizes the grid dimensions $(\Ncoarse,\Nfine,\Npsi) = (1024, 4096, 4096)$, and was performed with our CUDA graphics card implementation. The result shows, that we are able to go beyond classical limits given by the recurrence time and filamentation in plasma does not pose a major challenge to our methodology.

\begin{table}[htp!]
    \centering
    \begin{tabular}{l l}
    \toprule
        name                    & value \\
        \midrule
       $(\Ncoarse,\Nfine,\Npsi)$    & $(1024, 4096, 4096)$\\
        \texttt{inc. threshold} $\delta_\text{det}$     & $5\cdot 10^{-3}$ \\
        \texttt{stencil distance} of GALS  & $5 \cdot 10^{-3}$ \\
        time step size \rev{(CFL)} $\Delta t$                  & $1/(\Npsi L_v)$\\  
        \bottomrule
    \end{tabular}
    \caption{Parameters of the CMM for the fine-scale computations.}
    \label{tab:CMM-fine-scale-computations}
\end{table}

\subsubsection{Nonlinear Landau Damping}
Next, we present simulations of nonlinear Landau damping and compare the results to \cite{ParkerDellar2015,heath2012discontinuous}. For this simulation, we use the same CMM settings listed in \cref{tab:CMM-fine-scale-computations}, as they have shown to be free of recurrence time aliasing effects. The solution is initialized as in \cref{eq:landau_damping_inicond} and all the physical parameters are set as listed in \cref{tabl:landaudamping}.

\begin{figure}[htp!]
    \centering
    \begin{subfigure}{0.48\linewidth}
    \setlength\figureheight{0.6\linewidth}%
    \setlength\figurewidth{0.6\linewidth}%
    \caption{$k=0.5$}
    \includetikz{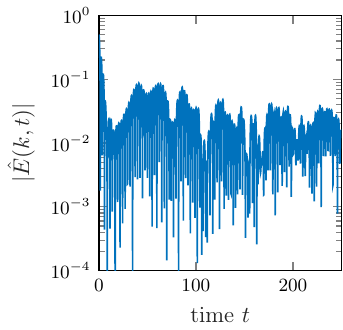}
    \label{fig:non-ld_Emodes1}
    \end{subfigure}
    \begin{subfigure}{0.48\linewidth}
    \setlength\figureheight{0.6\linewidth}%
    \setlength\figurewidth{0.6\linewidth}%
    \caption{$k=1$}
    \includetikz{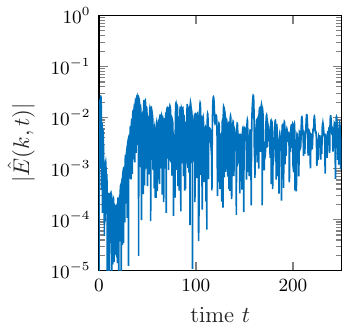}
    \label{fig:non-ld_Emodes2}
    \end{subfigure}
    
    \begin{subfigure}{0.48\linewidth}
    \setlength\figureheight{0.6\linewidth}%
    \setlength\figurewidth{0.6\linewidth}%
    \caption{$k=1.5$}
    \includetikz{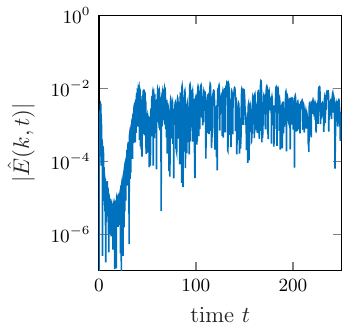}
    \label{fig:non-ld_Emodes3}
    \end{subfigure}
    \begin{subfigure}{0.48\linewidth}
    \setlength\figureheight{0.6\linewidth}%
    \setlength\figurewidth{0.6\linewidth}%
    \caption{$k=2$}
    \includetikz{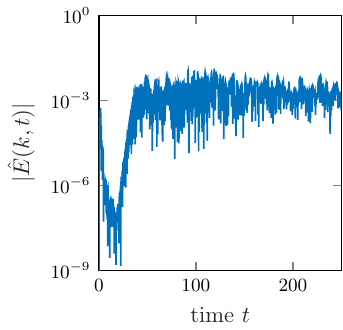}
    \label{fig:non-ld_Emodes4}
    \end{subfigure}
    
    \caption{Amplitude $\abs{\hat{E}(k,t)}$ of the first four modes of the electric field of the non-linear Landau damping.}
    \label{fig:non-lin-4fmodes}
\end{figure}

\begin{figure}[htp!]
    \centering
    \setlength\figureheight{0.8\linewidth}%
    \setlength\figurewidth{0.8\linewidth}%
    \includetikz{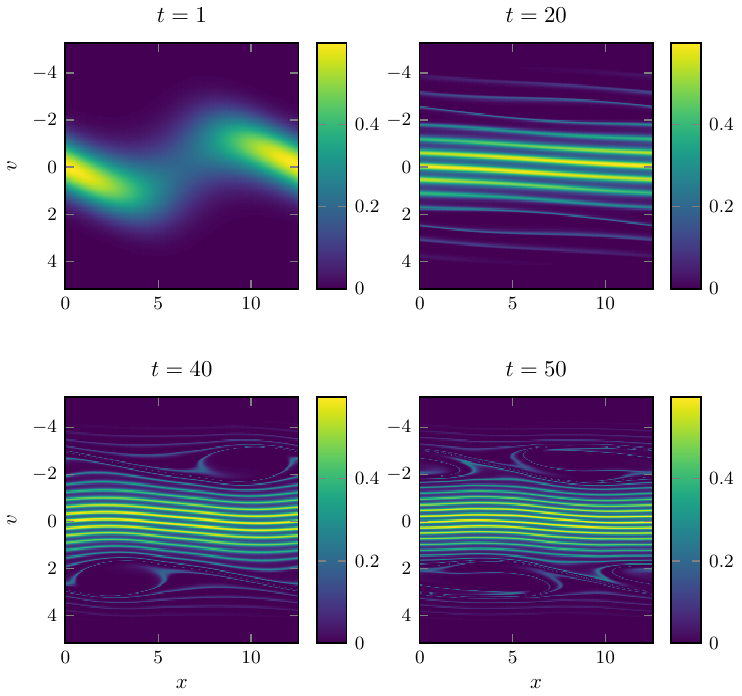}
    \caption{PDF of the non-linear Landau damping for four different time instances. Fine-scale computation with parameters listed in \cref{tab:CMM-fine-scale-computations}.}
    \label{fig:landau_damping_PDF_non_lin}
\end{figure}

As the perturbation is increased to $\epsilon=0.5$ the higher Fourier modes in the electric field get excited. We show the first four Fourier modes that correspond to the wave numbers $k=0.5,1,1.5,2$ in \cref{fig:non-lin-4fmodes} sampled till $T=250$. Confirming the observations of \cite{heath2012discontinuous,ParkerDellar2015} we see decreasing maximal amplitudes of the Fourier modes approximately until $t\approx 15$. Thereafter the amplitudes of all four modes increase until they reach their maxima at $t\approx 40$. The PDF is shown for four different time instances in \cref{fig:landau_damping_PDF_non_lin}.
We note that by closer examination the general behavior of the four modes is closer to the simulations shown in \cite{ParkerDellar2015}.
From the fits in \cref{fig:landau_damping_epot_non_lin} we can determine the two characteristic damping values $\gamma_1=-0.2872$ and $\gamma_2=0.0865$ that are in agreement with the values in the literature $\gamma_1=-0.287$ estimated by
\cite{heath2012discontinuous}, $(\gamma_1,\gamma_2)=(-0.281,0.084)$ estimated by \cite{ChenKnorr1976}.

\subsection{Two-Stream Instability}

The two-stream instability with an initial condition consisting of two counter-propagating cold electron beams initially located at $v_0 =3$
is considered in the following.

The initial condition is given by
\begin{equation}
    f_0(x,v) = (1+\epsilon\cos(kx))\frac{1}{2\sqrt{2\pi}} (e^{-(v-v_0)^2/2}+e^{-(v+v_0)^2/2})
\end{equation}
and the parameters are listed in \cref{tabl:2stream}.

\begin{table}[htp!]
    \centering
    \begin{tabular}{l l}
    \toprule
       name  & value \\
       \midrule
                \rev{spatial domain}    &     \rev{$\domain_x\in[0, 2\pi/k]$}\\
         \rev{velocity domain} & \rev{$\domain_v\in[-5\pi,5\pi]$}\\
         \rev{effective vel. domain} & \rev{$\domain_{v^*}\in[-4.75\pi,4.75\pi]$}\\
         wave number $k$ &    $0.2$\\
         $\epsilon$     &   $5 \times 10^{-2}$ \\ 
         $v_0$           &       $3$\\
         $T$              &     80\\
         \bottomrule
    \end{tabular}
    \caption{Parameters taken from \cite{PhamHelluyCrestetto2013} for the two-stream instability.}
    \label{tabl:2stream}
\end{table}

\paragraph{Convergence Tests.} We first perform basic convergence tests in space and time, similar to the ones in \cref{subsec:linear_landau_damping}.
For the convergence in time, we vary the spatial grid size of the map, while keeping all other maps at the same resolution. After $t=10$ we stop the simulation and compare it to the reference solution. Convergence in time is investigated by comparing a reference solution at $\Delta t=1/256$ to consecutive decreasing time step sizes. The different grid and time step sizes for our studies are listed in \cref{tabl:convergence-params}. The results of the spatial convergence shown in \cref{table:EOC-2stream} as well as the temporal convergence \cref{table:two_stream_time_conv} indicate that the method is of 3rd order in space and time.

\begin{table}                                                             
\centering               
\begin{tabular}{cccccc}   
\toprule                 
$\Delta t$ & $\Delta f$ &$\Delta \mathbf{X}$ \\
\midrule                 
$L/16$ & - & -  \\       
$L/32$ & 2.8 & 3.0  \\   
$L/64$ & 2.9 & 3.0  \\   
$L/128$ & 2.8 & 2.7 \\   
\bottomrule              
\end{tabular}          
\caption{Temporal order of convergence EOC as defined in \cref{eq-def:EOC} for the two-stream test case. The reference is CMM with $L/256$.}          
\label{table:two_stream_time_conv}
\end{table}

For a sanity check, we also compare the results to the spectral method introduced in \cref{sec:pseudo-spectral}. We see 3rd order convergence when comparing the CMM computations to the reference solution of the spectral method that is computed on a $512 \times 512$ grid. 

\begin{table}                              
\centering                 
\begin{tabular}{ccccccc}   
\toprule                   
$1/N$ & $\Delta f$ & $\Delta \mathbf{X}$ & $\Delta M$ & $\Delta P$ & $\Delta \Etot$ & $\norm{f-f^\text{spect}}_\infty$ \\
\midrule                   
1/32 & -    & -   & -   & -  & - & -  \\                                 
1/64 & 3.2 & 1.8 & 2.9 & 2.9 & 2.9 & 3.1\\                               
1/128 & 2.9 & 2.0 & 2.9 & 2.9 & 2.9 & 2.8 \\                             
1/256 & 3.1 & 2.4 & 3.1 & 3.1 & 3.0 & 2.9\\                              
1/512 & 3.4 & 2.9 & 3.3 & 3.3 & 2.9 & 3.0\\                              
\bottomrule                                                              
\end{tabular}   
\caption{Experimental order of convergence EOC as defined in \cref{eq-def:EOC} for the two-stream instability.}        
\label{table:EOC-2stream}                 
\end{table}

\paragraph{Fine-scale Structures and Zoom Properties.}

\begin{figure}[htp!]
    \centering
    \setlength\figureheight{0.7\linewidth}%
    \setlength\figurewidth{1\linewidth}%
    \includetikz{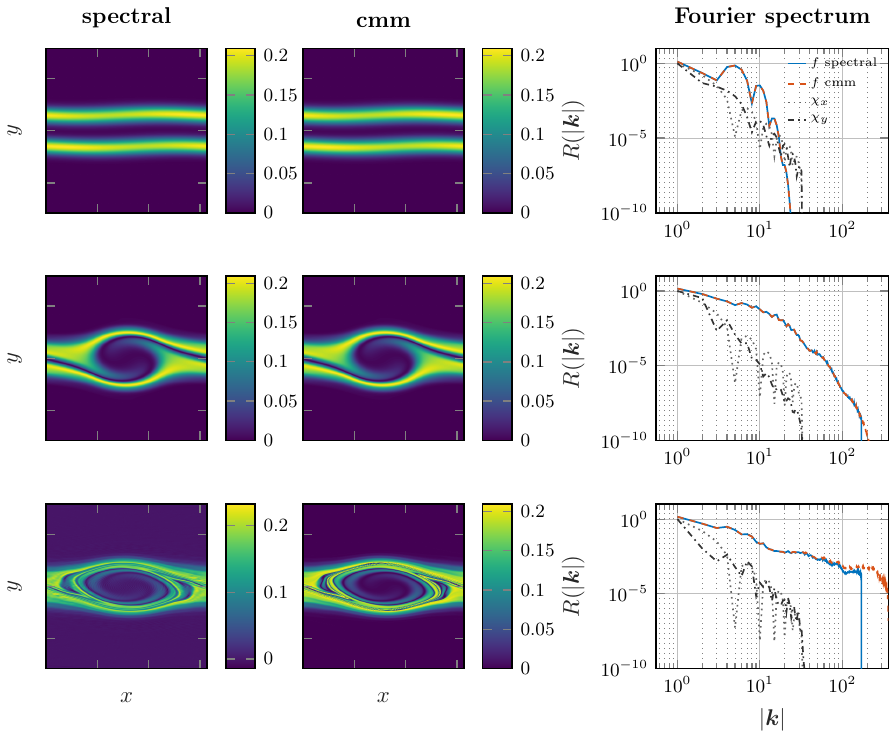}
    \caption{Time evolution of $f$ for the two-stream instability, at $t=1, 15, 40$. Resolution: spectral 512, CMM $\Ncoarse=64,\Nfine=\Npsi=512$,
    right: Radial Fourier spectra $R(\abs{\vect{k})}) = \int_{\abs{\vect{k}}-\Delta k}^{\abs{\vect{k}}+\Delta k}\abs{\mathcal{F}[g]({\vect{k})}}\d \vect{k}$, where $g$ is either the snapshot of the current submap $\chi_{x/y}$ in $x/y$ direction or the PDF.}
    \label{fig:2stream}
\end{figure}

One 
advantage of the CMM
method compared to the pseudo-spectral and other grid-based methods is that fine scales are not advected directly, but carried throughout the map. Thus the fine scales do not need to be
resolved during the advection, allowing for a coarser map. This is illustrated in \cref{fig:2stream}. The figure shows a time evolution of the PDF for the two-stream instability. It directly compares the pseudo-spectral to the CMM computation in the physical and in the Fourier domain. 
\rev{Over time, small-scale patterns form, due to plasma filamentation, filling up the Fourier spectra. At $t=40$, all Fourier modes are filled, causing spurious oscillations in the spectral simulation.}
Instead, CMM allows for an accurate representation even so the fine grid has a full Fourier spectrum because the grid of the sub-maps remains well-resolved as indicated by the fast decay of the Fourier spectra, shown on the right column of \cref{fig:2stream}.

\begin{figure}[htp!]
    \centering
    \includegraphics[width=1\textwidth,trim=0 0 50 0,clip]{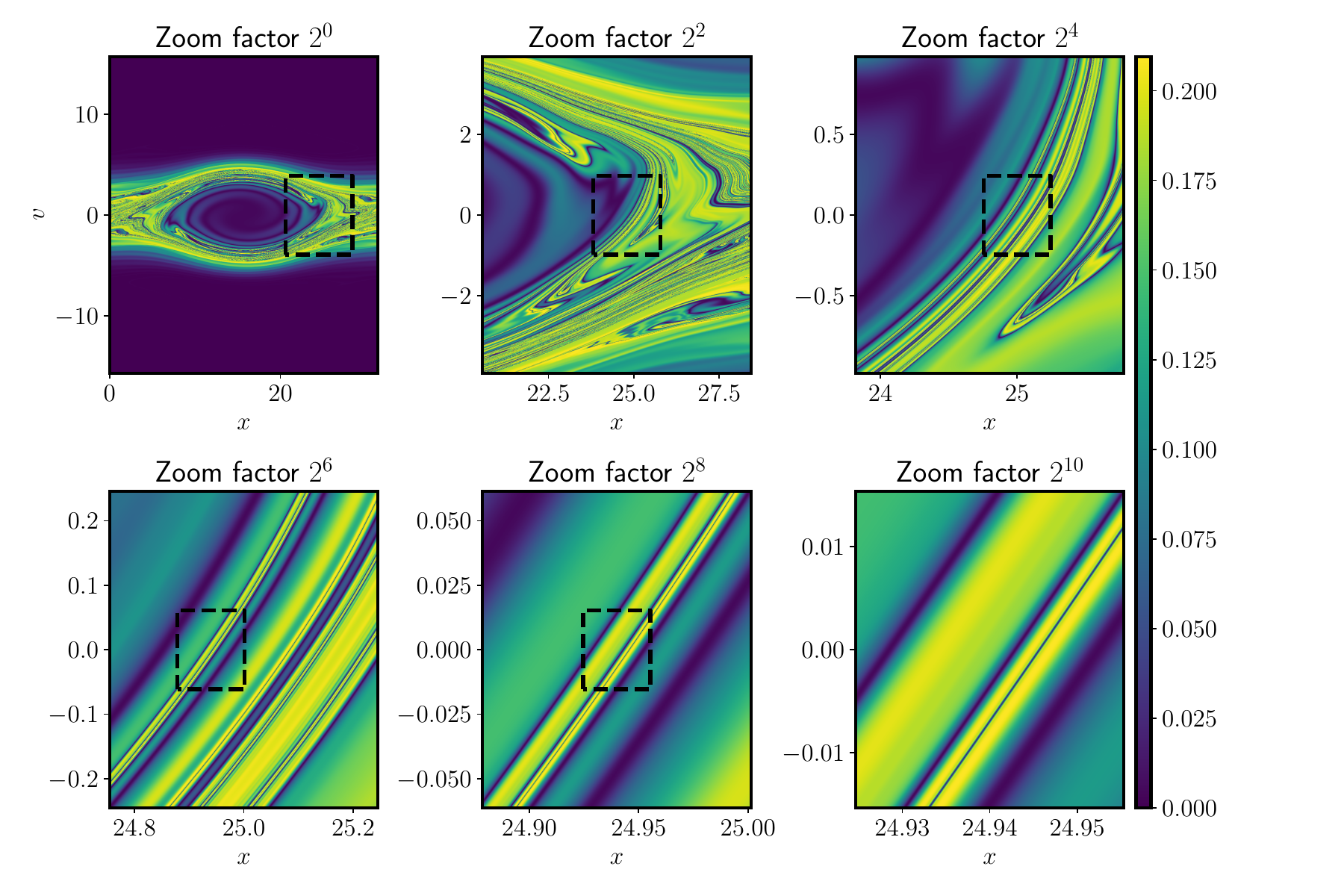}
    \caption{Two-stream instability computed with CMM, $\Ncoarse=1024,\Nfine=\Npsi=4096$. Successive zooms of the distribution function $f$ at $t=80$ to illustrate the multitude of resolved scales. The region inside the dotted black square is successively amplified by a factor 4 in each direction up to a factor of $2^{10}$.}
    \label{fig:2stream_zoom}
\end{figure}

To showcase the capability of the method we simulated the two-stream instability on a $(\Ncoarse,\Nfine,\Npsi) = (1024, 4096, 4096)$ configuration up to time $t=80$ and performed a zoom on the resulting solution, by interpolating the map at the corresponding zoom window and evaluating the initial condition on this window. \Cref{fig:2stream_zoom} shows the ability to resolve even the finest structures, using a zoom of factor $2^{10}$ compared to the initial domain size.

\section{Conclusion}
\label{sec:conclusion}


We developed the characteristic mapping method for the one+one dimensional Vlasov--Poisson equation. 
The method is an extension of the previously developed algorithm for the incompressible Euler equations \citep{yin2021characteristic, yin2023characteristic}. Differences exist concerning the transport term, which also contains divergence-free velocities in phase space, but these are no longer periodic in all directions. Thus we have used a periodization trick to be able to apply the CMM directly. 
Numerical analysis of the method showed third-order convergence in space and time. 
Numerical experiments confirmed our theoretical results for different test cases, including linear and non-linear Landau damping and the two-stream instability.
We performed extensive benchmarking and compared the CMM results with those published in the pertinent literature and moreover with a Fourier pseudo-spectral method. Thus the validity of CMM has been shown and its efficiency has been assessed. 

The semi-Lagrangian structure of CMM plays a pivotal role in addressing the issue of filamentation, particularly with regard to the temporal discretization induced by the CFL criteria. By carefully adjusting the sampling grid resolution, we have demonstrated that CMM can effectively eliminate filamentation problems. Notably, our investigations have revealed that the method is capable of resolving Landau damping to machine precision, all while remaining immune to recurrence time effects.
In comparison to other semi-Lagrangian methods, CMM stands out for its inherent bounds preservation. Unlike competing methods that often require additional measures such as positivity limiters or filters, CMM excels in preserving the integrity of physical bounds without the need for such supplements \cite{ZerroukatWoodStaniforth2005, ZerroukatWoodStaniforth2006,CrouseillesMehrenbergerSonnendruecer2010}.
CMM's compositional structure is another remarkable aspect of the method. It enables us to resolve fine-scale structures at varying resolutions through the application of zooms, a unique feature not shared by any other semi-Lagrangian approach. However, it is worth noting that this advantage comes with the requirement of storing sub-maps and executing compositions during runtime.
Our implementation of CMM also demonstrates computational efficiency, yielding time savings \rev{of one order of magnitude} when compared to a pseudo-spectral implementation of the method. These findings underscore the practical benefits of adopting CMM for numerical simulations.
Furthermore, our research indicates that CMM exhibits excellent conservation properties with third-order convergence. While the method showcases strong conservation characteristics, full numerical conservation of mass, energy, and momentum is a subject of future work.
Furthermore, future work will address VP using CMM in higher dimensions.

\section*{Author Contribution Statement (CRediT)}


In the following, we declare the author's contributions to the publication:

\vspace{5pt}
{\small
\noindent
\begin{tabular}{@{}lp{10.8cm}}
\textbf{Philipp Krah:} & implementation MATLAB \& C\texttt{++}/CUDA code, writing initial draft, numerical studies, visualization \\
\textbf{Xi-Yuan Yin:} &  implementation MATLAB, numerical analysis, boundary periodification, writing initial draft, reviewing \& editing\\
\textbf{Julius Bergmann:} & implementation C\texttt{++}/CUDA code, reviewing \& editing\\
\textbf{Jean-Christophe Nave:} & initial idea, editing and supervision\\
\textbf{Kai Schneider:} & initial idea, writing the initial draft, editing and supervision, funding acquisition, project administration \\
\end{tabular}
}

\section*{Conflict of Interest}
The authors declare that they have no conflict of interest.

\section*{Code and Data Availability}
All scripts to reproduce the results are available at:

\begin{centering}
\url{https://github.com/orgs/CharacteristicMappingMethod/repositories}.
\end{centering}

\section*{Acknowledgement}
The authors were granted access to the HPC resources of IDRIS under allocation No. AD012A01664R1 attributed by Grand Équipement National de Calcul Intensif (GENCI).
Centre de Calcul Intensif d’Aix-Marseille is acknowledged for granting access to its high-performance computing resources.
The authors acknowledge partial funding from the Agence Nationale de la Recherche (ANR), project CM2E, grant ANR-20-CE46-0010-01.
\rev{KS acknowledges financial support from the French Federation for Magnetic Fusion Studies (FR-FCM) and the Eurofusion consortium, funded by the Euratom Research and Training Programme under Grant Agreement No.633053. The views and opinions expressed herein do not necessarily reflect those of the European Commission.}

\bibliographystyle{abbrv}
\bibliography{references}

\newpage
\section*{Appendix}
\appendix

\section{Periodic Velocity Extension}
\label{appdx:BLayer}

We want to define a periodic velocity field $u_1(v)$ on the velocity interval $v\in [-L_v, L_v]:=\domain_v$ such that on the support of $f$, with $\text{supp}(f) \subset \domain_{v^*}:=(-v^*, v^*) \subset \domain_v$, we have $u_1(v) = v$. It follows that the derivative $\partial_v u_1$, the $u_1$ contribution to vorticity must be an approximation of the box function 
\begin{equation}
    B(v) = 1 - \frac{L_v}{L_v-v^*}\mathbf{1}_{|v| \in [v^*, L_v]},
\end{equation}
such that the integral of $B$ vanishes in order to satisfy the periodicity condition on $u_1$. In doing so, we add a padding region $[v^*, L_v]$ corresponding to a horizontal jet in order to make $u_1$ vanish smoothly at $v=L_v$. The length of the extension, $a = L_v-v^*$, can be seen as the amount of wasted computational resource and should be minimized subject to error tolerances. The function $B$ is only a generic example of a periodic extension, smoother extensions exist with the only requirement being that $B(v) \equiv 1$ in $\domain_{v^*}$ and the integral of $B$ is 0. We note that due to the requirement that $B$ is constant in an open interval, analytic functions are out of question and the best extensions we can hope for are in the $C^\infty$ class. The regularity of the extension will affect the decay of the coefficients of the Fourier transform and thereby the numerical error introduced. The general optimization problem is as follows: Given a family of extension $B(v; a)$ parametrized by $a$ ($v^*$ given, thereby also parametrized by $L_v$), for each computational grid size $N$, we minimize the total numerical error
\begin{equation}
    \mathcal{E}(N, a) = \mathcal{D}\left( \frac{N}{v^*+a} \right) + \mathcal{C}(a, N),
\end{equation}
where $\mathcal{D}$ is the spatial error of the numerical scheme which depends on the Nyquist frequency $N/L$ and $\mathcal{C}$ is the numerical error introduced by the velocity modification which we assume to be the on the order of the $k^{-1} \hat{B}(N; a)$, the size of the first truncated term in terms of the Fourier transform. We note that this is an over-estimation of $\mathcal{C}(a, N)$, indeed, assuming $B \equiv 1$ in $\domain_{v^*}$ the extension only introduces error outside the support of $f$ which then propagates towards the support itself through successive time stepping.

\begin{figure}[htp!]
    \centering
    \begin{subfigure}{0.48\linewidth}
    \centering
     \setlength\figureheight{0.6\linewidth}%
    \setlength\figurewidth{0.8\linewidth}%
    \includetikz{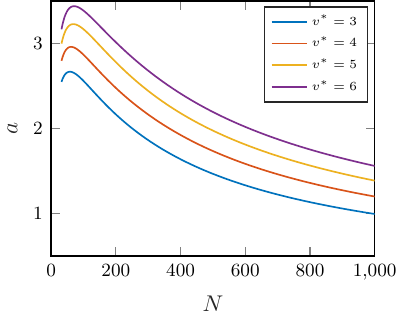}
    \caption{Length of the extension}
    \end{subfigure}
    \begin{subfigure}{0.48\linewidth}
    \centering
     \setlength\figureheight{0.6\linewidth}%
    \setlength\figurewidth{0.8\linewidth}%
    \includetikz{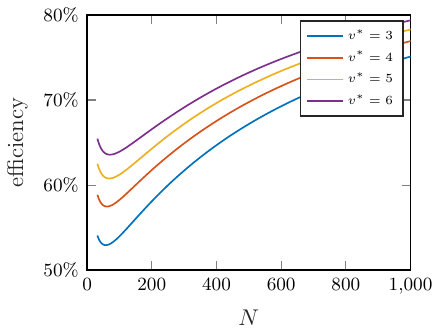}
    \caption{Efficiency of the extension}
    \end{subfigure}
    \caption{Choice of extension region optimized according to \cref{eq:ExtensionOptimizationFormula}. The efficiency is defined as the ratio of the number of grid points in the support against the total number of grid points.}
    \label{fig:extension_efficiency}
\end{figure}

We consider the simple extension using a bump function:
\begin{equation}
    B(v; a) = 1 - \frac{L}{a^2} \eta \left( \frac{\abs{x}-L_v}{a} \right) ,
\end{equation}
where $\eta$ is the standard mollifier $\eta(x) = \exp \left( \frac{-1}{1-x^2} \right)/I$ which satisfies $\text{supp}(\eta) = (-1, 1)$ and $\int_{-1}^1 \eta \d x = 1$. The technical note \cite{johnson2015saddle}, gave a very accurate asymptotic decay rate of $L_v a^{-7/4} k^{-3/4} \exp(-\sqrt{ak} ) $ for these mollifiers, from which we derive an estimate for the error introduced by the extension:
\begin{equation}
    \mathcal{C}(a,N) \propto L_v \left( \frac{aN}{v^*+a} \right)^{-7/4} \exp \left( - \sqrt{ \frac{aN}{v^*+a}} \right) .
\end{equation}
For a Hermite cubic implementation of the characteristic mapping method, the local spatial truncation error is $4^{th}$ order, i.e. assuming $\mathcal{D} \propto \left( \frac{N}{a+v^*} \right)^{-4}$, we \rev{obtain} the following approximate optimality condition for $a$ for a given fixed $N$:
\begin{equation} \label{eq:ExtensionOptimizationFormula}
    0 = \frac{4 (a + v^*)^3}{N^4} - ( N_a )^{-\frac{7}{4}} \exp \left( -  \sqrt{N_a} \right) \left( \frac{v^*}{4a} \left( 2 \sqrt{N_a} + 7  \right)   - 1  \right),
\end{equation}
where $N_a = \frac{aN}{a + v^*}$ is the number of wasted grid points. We solve this numerically and show the size of the optimized extension and its computational efficiency $\frac{v^*}{a + v^*}$ in \cref{fig:extension_efficiency}.

\section{Statistics of Fine-Scale Computations}
\label{appx:fine-scale-stats}

\begin{figure}[hbp!]
    \centering
    \setlength\figureheight{1.25\linewidth}%
    \setlength\figurewidth{0.9\linewidth}%
    \includetikz{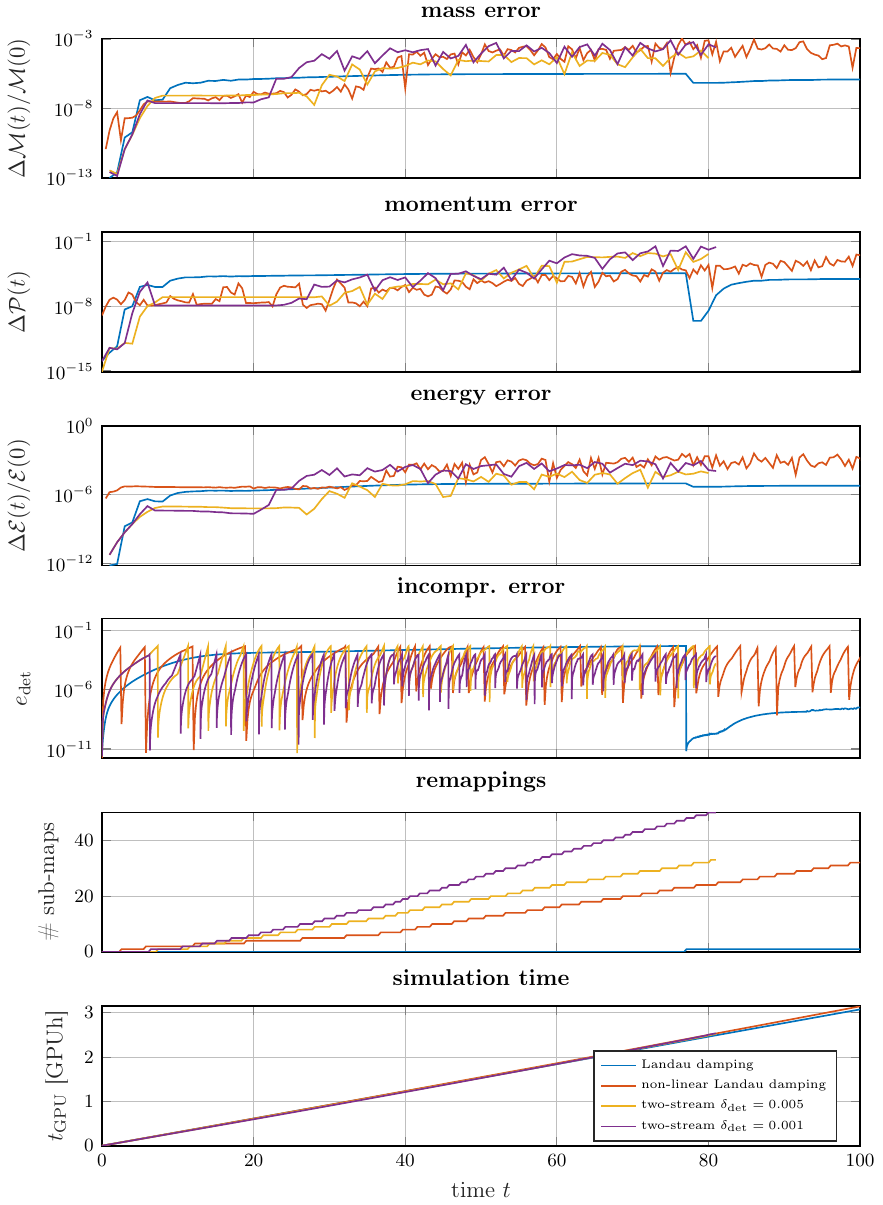}
    \caption{Statistical quantities of the fine resolved CMM computations: Landau damping, non-linear Landau damping and the two-stream instability. The legend is given in the lowest figure. The simulation parameters are listed in \cref{tab:CMM-fine-scale-computations}. Shown are (from top to bottom) relative mass error, momentum error, relative energy error, incompressibility error of the map, number of sub-maps, GPU time in hours on a \texttt{NVIDIA A100 SXM4 80 GB GPU.}
    \label{fig:CMM-statistics}}
\end{figure}

For the linear Landau damping, non-linear Landau damping and two-stream instability test case, we have performed fine scale CMM computations with algorithmic parameters listed in \cref{tab:CMM-fine-scale-computations} and physical parameters in \cref{tabl:landaudamping,tabl:2stream}. Results have been shown for the linear and non-linear Landau damping in \cref{fig:Landau_damping_without_recurence,fig:landau_damping_PDF_non_lin,fig:non-lin-4fmodes} and the zooms of the two-stream instability in \cref{fig:2stream_zoom}.
In \cref{fig:CMM-statistics} we show the time history of the conserved quantities $\mass,\Etot,\momentum$ and additional quantities like the compressibility error, the number of sub-maps utilized during the computations and the overall GPU-time.
For all fine-scale simulations, we use the same incompressibility threshold $\delta_\mathrm{det}=0.005$. However, we include one fine-scale simulation of the two-stream instability using $\delta_\mathrm{det}=0.001$ to visualize the difference in performance.
\Cref{fig:CMM-statistics} shows that all conserved quantities grow until a saturated state. Remappings affect the simulation results as can be seen from the direct comparison of the two-stream instability computations. After the first remapping of the simulation with $\delta_\mathrm{det}=0.001$ at  $t\approx 7$, the results compared to the simulation with $\delta_\mathrm{det}=0.005$ differ. However, the higher number of remappings in the case of $\delta_\mathrm{det}=0.001$ does not imply a significant rise in GPU time.
Nevertheless, it leads to higher consumption in memory, which is the limiting factor when choosing small values of $\delta_\mathrm{det}$.

\begin{figure}[hbp!]
    \centering
    \setlength\figureheight{1.25\linewidth}%
    \setlength\figurewidth{0.9\linewidth}%
    \includetikz{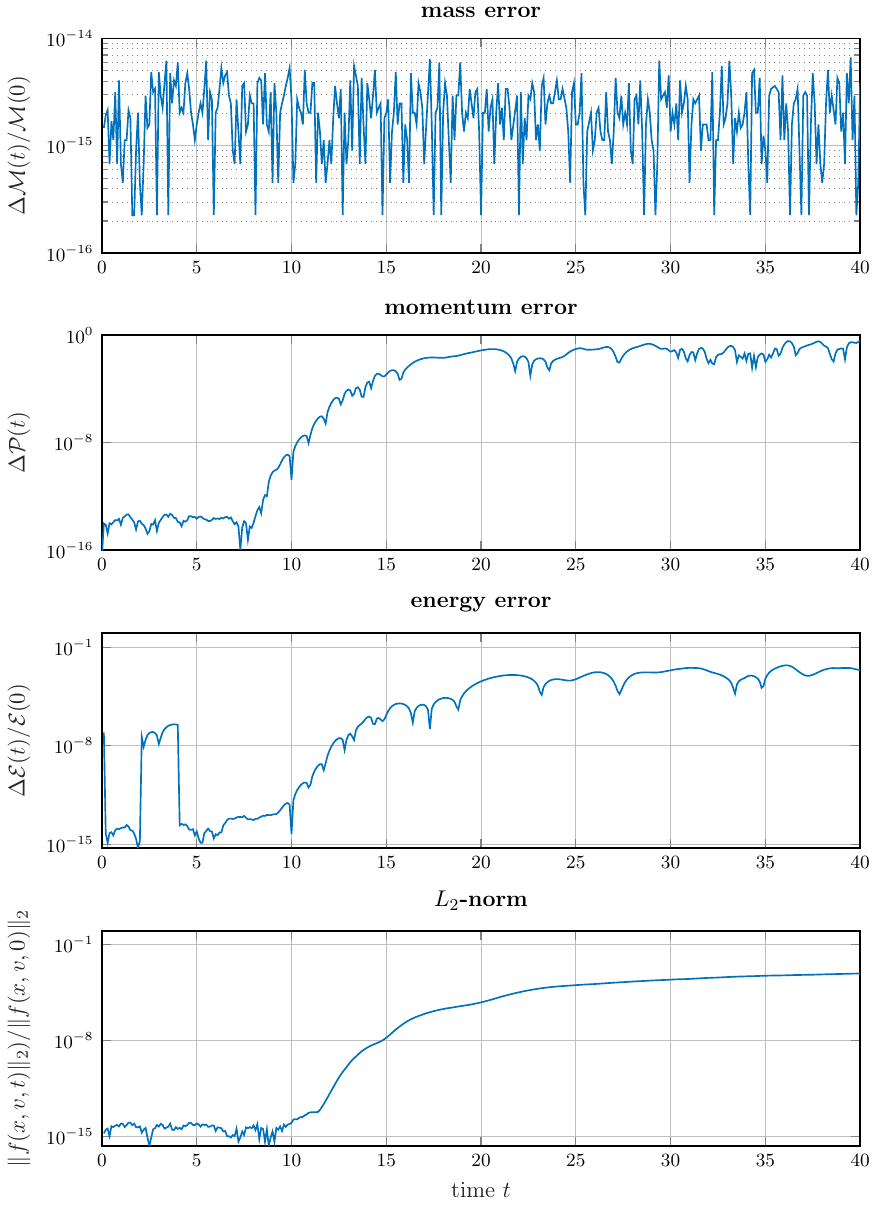}
    \caption{Statistics of spectral computation of the two-stream instability at $N\times N =512 \times 512$.  Shown are (from top to bottom) relative mass error, momentum error, relative energy error and $L_2$-error \rev{of $f$}.}
    \label{fig:spectral-statistics}
\end{figure}

Additionally to the CMM computations, we show the statistics of the pseudo-spectral computations in \cref{fig:spectral-statistics}. These are the statistics of the simulation shown in \cref{fig:2stream}.
Note that the total mass is conserved up to machine precision. $L_2$-norm of $f$ is preserved until the first oscillations in the simulation occur.

\end{document}